\RequirePackage[l2tabu, orthodox]{nag}
\documentclass[a4paper,english,oneside]{scrartcl}
\pdfoutput=1
\usepackage[english]{babel}
\usepackage[utf8x]{inputenc}
\usepackage{amsmath,amssymb,amsthm}
\usepackage{mathrsfs}
\usepackage[all]{xy}
\usepackage[binary-units=true]{siunitx}
\usepackage[square,sort,comma,numbers]{natbib}		
\usepackage{graphicx}
\usepackage{algpseudocode}
\usepackage{algorithm}
\usepackage{subfig}
\usepackage{mdframed}
\usepackage{geometry}
\usepackage{epstopdf}
\graphicspath{{plotting/}}

\newcommand{\N}{\mathbb{N}}
\newcommand{\Z}{\mathbb{Z}}
\newcommand{\Q}{\mathbb{Q}}
\newcommand{\R}{\mathbb{R}}
\newcommand{\C}{\mathbb{C}}
\newcommand{\rd}[1]{\lfloor #1 \rceil}			
\newcommand{\PP}{\R[\underline X]}			
\newcommand{\SOS}{\Sigma}				
\newcommand{\sls}{\mathcal{L}_{\leq}}			
\newcommand{\subproblems}{\mathcal{L}}			
\newcommand{\sphere}{\mathbb{S}}			
\newcommand{\cone}{\mathscr{C}}				
\newcommand{\stacktext}[1]{\stackrel{\text{#1}}}        
\newcommand{\uniDist}{\mathcal{U}}			
\newcommand{\NULL}{\texttt{NULL}}

\newcommand{\define}[1]{\textit{#1}}			
\newcommand{\stress}[1]{\textit{#1}}			

\newcounter{thmcounter}
\theoremstyle{plain}
\newtheorem{prp}[thmcounter]{Proposition}
\newtheorem{thm}[thmcounter]{Theorem}
\newtheorem{lem}[thmcounter]{Lemma}
\newtheorem{cor}[thmcounter]{Corollary}
\newtheorem{obs}[thmcounter]{Observation}
\theoremstyle{definition}

\newtheorem{rem}[thmcounter]{Remark}
\newtheorem{nota}[thmcounter]{Notation}

\algnewcommand\algorithmicinput{\textbf{input}}
\algnewcommand\Input[1]{\State \algorithmicinput\ #1}
\algnewcommand\algorithmicoutput{\textbf{output}}
\algnewcommand\Output[1]{\State \algorithmicoutput\ #1}
\algnewcommand\algorithmicprint{\textbf{print}}
\algnewcommand\Print[1]{\State \algorithmicprint\ \texttt{#1}}
\renewcommand{\algorithmiccomment}[1]{\bgroup\hfill\textit{//~#1}\egroup}

\title{Norm Bounds and Underestimators for Unconstrained Polynomial Integer Minimization}

\author{Sönke Behrends\footnote{Institute for
Numerical and Applied Mathematics, University of G\"ottingen, Lotzestr. 16--18, 37083 G\"ottingen,
Germany} 
\and Ruth Hübner\footnotemark[1] \and Anita Schöbel\footnotemark[1]}

\date{January 19, 2015}
\begin{document}
\maketitle

\begin{abstract}
We consider the problem of minimizing a polynomial function over the integer lattice. 
Though impossible in general, we use a known sufficient condition for the existence of 
continuous minimizers to guarantee the existence of integer minimizers as well. 
In case this condition holds, we use sos programming to compute the radius of a $p$-norm
ball which contains all integer minimizers. We prove that this radius is smaller than the radius
known from the literature.  Furthermore, we derive a new class of underestimators of the polynomial function. Using a Stellensatz from real algebraic geometry and again sos programming, we optimize over this class to get a strong lower bound on the integer minimum.

Our radius and lower bounds are evaluated experimentally. They show 
a good performance, in particular within a branch and bound framework.
\end{abstract}

\paragraph*{Keywords} 
integer optimization; polynomials; lower bounds; branch and bound

\section{Introduction}
\label{sec:intro}
Given a multivariate polynomial $f:\R^n\rightarrow \R$, we consider its minimization over the integer lattice, i.e., the problem
\begin{equation}
\begin{aligned}
\min~ & f(x) \\ 
\text{s.t.}~& x\in \Z^n.
\end{aligned}
\label{eq:problem} \tag{IP}
\end{equation}
This is a special type of a nonlinear integer optimization problem and is incomputable in general: Hilbert's tenth problem asks if there exists an algorithm that decides whether for a given polynomial $f$ with integer coefficients the equation $f(x) = 0$ has a solution $x~\in~\Z^n$. Seventy years later it was proved by Matiyasevich~\cite{matiyasevich1970enumerable} that no such algorithm can exist. So if there was an algorithm to solve \ref{eq:problem}, we would also get an algorithm to decide whether $f(x)=0$ has an integer solution by minimizing $f^2$ over $\Z^n$. Consequently, \ref{eq:problem} cannot be solved for general polynomials $f$. In this paper 
we consider a subclass that leads to solvable problems.

\subsection{Outline}
Once the notation, a little background on sos (sum of squares) programming and a Stellensatz from algebraic geometry are introduced (Section~\ref{sec:preliminaries}), we review, in order to make the problem tractable, a sufficient criterion from the literature for the existence of \stress{continuous} minimizers (Section~\ref{sec:existence_of_minimizers}). This criterion actually holds for integer minimzers, too: Integer minimizers exist if the highest order terms of $f$ attain positive values on $\R^n~\setminus~\{0\}$, we say that the \define{leading form} of $f$ is \define{positive definite}. But deciding positive definiteness is NP hard, hence we approximate this problem by sos programming. However, this only tells us that minimizers exists, but not where they are located. We locate the minimizers by computing -- again using sos programming -- the radius of a $p$-norm ball that contains all integer minimizers (Section~\ref{sec:norm_bound}), or simply \define{norm bounds} on the minimizers. In principle, once a norm bound is known, \ref{eq:problem} is solvable by enumeration.
We proceed by deriving a class of polynomials with obvious integer minimizer (Section~\ref{sec:underestimators}) that serve as underestimators to $f$. Using sos programming, we may choose the underestimator $g$ with the strongest lower bound. Firstly, we search for a global underestimator which is later refined to underestimation on sublevel sets, yielding stronger bounds. This refinement further allows to prove that, provided $f$ has a positive definite leading form, there always are underestimators in our class that can be found by sos programming.
To find the optimal solution to \ref{eq:problem}, instead of enumeration, we use an underestimator $g$ from our class to obtain lower bounds within a branch and bound approach (Section~\ref{sec:implementation}). We continue with an experimental evaluation of the norm bounds, of the lower bounds and of the performance of the underestimators within branch and bound on random instances. 
The paper ends with a conclusion and ideas for future research (Section~\ref{sec:ausblick}).

\subsection{Literature review}

The literature on nonlinear integer programming is vast. For an overview, a presentation of key techniques and complexity results as well as numerous references for further reading, see the article~\cite{hemmecke2010nonlinear}, which comes as chapter of~\cite{junger201050years}. For a recent survey on nonlinear \stress{mixed-integer} programming (a subset of the variables may be continuous), see~\cite{lee2012mixed}. 

\medskip

Throughout our work we rely heavily on methods from constrained \stress{continuous} polynomial optimization. Based on work of Shor~\cite{shor1987class,shor1997modified}, Parrilo~\cite{parrilo2000structured} suggested a method now known as sos programming that makes continuous polynomial optimization accessible to semidefinite programming (see, e.g.,~\cite{wolkowicz2000handbook} for the latter), whilst Lasserre~\cite{lasserre2001global} published the dual approach, based on moment sequences.
Since the emergence of the two ground-breaking publications by Parrilo and Lasserre, many results on continuous polynomial optimization via sos techniques and its theoretical background have been published: The expository paper~\cite{parrilo2003minimizing} shows that existing algebraic techniques are outperformed by the sos method. As in-depth treatments, we refer to~\cite{anjos2012handbook} for the interplay of semidefinite, conic and polynomial optimization, and~\cite{blekherman2013semidefinite} for a focus on the geometry involved.  For an algebraic treatment, we mention Marshall's book~\cite{marshall2008positive}. We point out Laurent's elegant survey~\cite{laurent2009sums}, which treats, among other aspects, the duality of the sos and moment approach.

\medskip

A special case of our problem, unconstrained quadratic integer minimization, is considered by \cite{Buch-Hueb-Sch13}. We did not find results in the literature that consider the unconstrained integer minimization problem for multivariate polynomials of arbitrary degree. 
Regarding nonlinear integer minimization with constraints or additional assumptions, integrality turns even seemingly simple problems incomputable: Using the aforementioned result of Matiyasevich, Jeroslow \cite{jeroslow1973there} proved that there cannot be an algorithm for integer minimization of a linear form subject to quadratic constraints.

\medskip

But substantial special cases are solvable, for example, every integer problem with a bounded feasibile set is solvable.
More specifically, an important case is boolean programming, see~\cite{boros2002pseudo} for a survey. A classic approach is linearization by introducing new variables and constraints (for early results see, e.g., \cite{fortet1960algebre}). In theory, also a general bounded integer polynomial optimization problem can be reduced to the binary case \cite{watters1967reduction}, but this is not practicable since the number of variables grows too quickly.
Another technique for boolean programming is the reduction to a quadratic problem which can be done with significantly fewer variables and constraints \cite{rosenberg1975reduction,buchheim2007efficient}.
Another substantial case that gained attention are (quasi-)convex problems, as the incomputability results do not hold for this case \cite{khachiyan1983convexity,khachiyan2000integer}. \cite{hildebrand2013new} present a Lenstra type algorithm for quasiconvex integer polynomial optimization. 

\medskip

For integer minimization of arbitrary polynomials, a common way of solving \ref{eq:problem} is branch and bound as proposed (originally only for convex functions) by~\cite{gupta1985branch}. A popular method is to calculate convex underestimators (see, e.g., \cite{lasserre2011convex}) to obtain lower bounds. As a different approach, if the feasible set is a box, \cite{buchheim2014box} compute separable underestimators wich give lower bounds that are easy to obtain. In contrast, \cite{de2006integer} directly compute lower and upper bounds, i.e., no underestimators, for nonnegative polynomials on polytopes.

\section{Preliminaries}
\label{sec:preliminaries}

In this section we collect some basic notation and facts as well as a few theorems from algebraic geometry that we use to derive our main results.

\subsection*{Notation and basic properties of polynomials}
\label{subsec:polynomials}

We write a polynomial $f:\R^n \to \R$ in $n$ unknowns $X_1, \ldots, X_n$ using multi-indices $\alpha = (\alpha_1, \ldots, \alpha_n) \in \N_0^{n}$ via $$f = \sum_{\alpha \in \N_0^n} a_\alpha X^\alpha = \sum_{\alpha \in \N_0^n}a_\alpha X_1^{\alpha_1} \cdots X_n^{\alpha_n},$$
for some unique coefficients $a_\alpha \in \R$, only finitely many nonzero, and monomials $X^\alpha := X^{\alpha_1} \cdots X^{\alpha_n}$.
The \define{modulus} of $\alpha \in \N_0^n$ is $| \alpha| = \alpha_1 + \ldots + \alpha_n$. With these conventions, the degree of $f$ is given by
$$\deg f := \max \left \{ | \alpha | \ \big | \  a_\alpha \neq 0 \right \}.$$
The \define{ring of polynomials} in the unknowns $X_1, \ldots, X_n$ is denoted by $\R[X_1, \ldots, X_n]$, which we 
abbreviate to $\PP$. We use $\underline X$ here in order to distinguish the multivariate from the univariate case.
\medskip

A polynomial $f$ is \define{homogeneous} if all monomials in $f$ have the same degree, say $d$. That is, $f$ is homogeneous if
$f = \sum_{| \alpha| = d} a_\alpha X^\alpha$. In this case one has
$$f(\lambda x) = \lambda^d f(x), \quad x \in \R^n,\ \lambda \in \R.$$
This implies that a homogeneous polynomial is uniquely determined by its values on any of the $p$-norm unit spheres
$$\sphere^{n-1}_p := \left\{ x \in \R^n \ \big | \ \| x \|_p = 1 \right \}, \quad p \in [1, \infty].$$

A homogeneous polynomial $f$ is \define{positive definite} if $f(x) > 0$ for $x \neq 0$. Similarly, a (possibly nonhomogeneous) polynomial $f$ is \define{positive semidefinite} if $f (x) \geq 0$ for all $x \in \R^n$, for short $f > 0$ and $f \geq 0$. For a homogeneous polynomial $f$ and some $p \in [1, \infty]$, we often use the following equivalent characterization:
\begin{equation}
\begin{aligned}
	\label{eq:definiteness_by_bounds}
	f \geq 0& \Longleftrightarrow  \exists c \geq 0:\ f(x) \geq c \text{ for all } x \in \sphere^{n-1}_p ,\\
	f > 0 	& \Longleftrightarrow  \exists c > 0: 	\ f(x) \geq c \text{ for all } x \in \sphere^{n-1}_p.
\end{aligned}
\end{equation}

A homogeneous polynomial is also called a \define{form}. Any polynomial $f \in \PP$ can be uniquely decomposed as
$$ f = \sum_{j=0}^d f_j$$	
where $d := \deg f$ and the $f_j$ are homogeneous polynomials of degree $j$, called the \define{homogeneous components} of $f$. 
The highest degree component, $f_d$, is called the \define{leading form} of $f$.
\medskip	

Given a vector $h \in \R^n$, we denote by $\rd{h}$ the vector resulting of rounding each component of $h$ to its nearest integer. Finally, we need the notion of a \define{sublevel set}: For a function $f: U \rightarrow \R$ from some set $U$, the \define{sublevel set} of \define{level} $z \in \R$ is defined by 
$$\sls^{f}(z) = \{x \in U \ | \ f(x) \leq z \}.$$

\subsection*{Nonnegativity and sums of squares}
Continuous minimization of a polynomial as well as deciding
non-negativity of a polynomial are well known to be NP-hard problems,
even if one fixes the degree to $d=4$
\cite{nesterov2000squared}. 
Deciding if $g$ is an underestimator of $f$ means to decide 
if $f-g$ is nonnegative. As this is NP-hard we
use a tractable sufficient criterion for nonnegativity in
the following:
We search for a
decomposition into a \define{sum of squares}, or \define{sos} for
short, which is a sufficient, but not necessary condition for
nonnegativity~\cite{marshall2008positive}. Formally, a polynomial $f
\in \PP$ is a sum of squares if there are $u_1, \ldots, u_l \in \PP$
such that $f = \sum_{i=1}^l u_i^2.$ We sometimes use the following
property:
\begin{lem}
\label{lem:sos_nonvanishing}
Suppose $v= u_1^2 + \cdots + u_k^2$
for some given $u_1, \ldots, u_k \in \R[\underline X]$ and $u_1 \neq 0$. Then $v \neq 0$, and 
$$\deg v = 2 \max_{1 \leq i \leq k} \deg u_i.$$
\end{lem}
\begin{proof}
See, e.g.,~\cite[Cor. 1.1.3]{marshall2008positive}.
\end{proof}

The convex cone 
$$\SOS := \left\{ f \in \PP \ \big |\ \exists u_1, \ldots, u_l \in \PP \text{ s.t. } f = \sum_{i=1}^l u_i^2 \right\}$$
in $\PP$ contains all polynomials $f$ which are sos in $\PP$. 
It is possible to optimize a linear form such that affine combinations of the decision variables and given polynomials lie in this cone: 
Such an \define{sos optimization problem} or \define{sos program} is tractable, as 
it is equivalent to a semidefinite program. For details, see e.g.~\cite{anjos2012handbook,blekherman2013semidefinite} on sos programming and~\cite{wolkowicz2000handbook} for semidefinite programming. Formally, an sos program has the form 

\begin{alignat}{2}\label{def:sosProgram}
\max	\quad 	& b_1y_1 + \cdots + b_m y_m \nonumber \\ 
\text{s.t.}	\quad 	& a_{i0} + y_1 a_{i1} + \cdots + y_m a_{im} \in \SOS,	\quad 	& & i = 1, \ldots, k, \\
		& y_i \in \R, \quad 						& & i = 1, \ldots, m, \nonumber
\end{alignat}
where $y_i \in \R$ are the decision variables, and $b_i \in \R$ 
as well as $a_{ij} \in \R[X_1, \ldots, X_n]$ are fixed. In our paper we
use sos programming for two purposes: 
To find an optimal underestimator, see Section~\ref{sec:underestimators}, and for constrained continuous minimization of polynomials as done at the end of Section~\ref{sec:preliminaries}. For both, we use a result from real algebraic geometry, known as Putinar's Stellensatz, outlined next. 

\subsection*{A result from algebraic geometry}
In this section we introduce Putinar's Stellensatz.
See~\cite{nie2007complexity} and the references therein for a discussion and the origins 
of the Stellensatz. 
For a finite collection of multivariate polynomials $S = \{g_1, \ldots, g_s \} \subset \PP$, define the semi-algebraic set $K_S$ as
\begin{equation}
K_S := \left \{ x \in \R^n \ | \ g_1(x) \geq 0, \ldots, g_s(x) \geq 0 \right \},
\label{eq:KS}
\end{equation}
where our notation follows~\cite{marshall2008positive}.
The Stellensatz we consider gives a sufficient conditions which allows to construct every polynomial $f \in \PP$ with $f > 0$ on $K_S$ from the given inequalities $g_i(x) \geq 0$. 
To this end, for $S$ as above, the \define{quadratic module} generated by $S$ is given by
\begin{equation}
M_S	:= \left \{ \sum_{i=0}^s \sigma_i g_i | \ \sigma_0, \ldots, \sigma_s \in \SOS \right \} 
\label{eq:syntacticQM}
\end{equation}
where $g_0:=1$. 
For the Positivstellensatz to hold we need $M_S$ to be \define{Archimedean}. This is the case if there is a polynomial $q \in M_S$ such that the set $K_{\{q\}} = \{ x \in \R^n \ |\ q(x) \geq 0 \}$ is compact. 

\begin{thm}[Putinar]
\label{thm:putinar_pss}
Let $M_S$ be Archimedean and $f \in \PP$. Then $f(x) > 0$ for all $x \in K_S$ implies $f \in M_S$.
\end{thm}

\subsection*{Lower bounds for constrained continuous minimization}

In Section~\ref{sec:existence_of_minimizers} we will see that in order to 
decide existence of minimizers, we need to compute a lower bound on the 
minimum of the leading form on the sphere $\sphere^{n-1}_p$. As the sphere 
is semi-algebraic 
for even $p$, sos methods can be applied to find such a lower bound. 

\medskip

In the following we describe how lower bounds on 
\begin{equation}
\begin{aligned}
\min	 	\quad & f(x) \\
\text{s.t.} 	\quad & x \in K_S, 
\end{aligned}
\label{eq:constrained_optimization}
\end{equation}

where $K_S$ defined as in (\ref{eq:KS}) can be derived by sos-programming. The method we outline follows Schweighofer~\cite{schweighofer2005optimization}, based on Lasserre's~\cite{lasserre2001global} work. 
We consider a hierarchy \ref{eq:constrained_optimization_sos}, $k = 1, 2, \ldots$, of sos programs 
\begin{equation}
\begin{aligned}
\label{eq:constrained_optimization_sos}
\max 	\quad & y_1 \\
\text{s.t.}	\quad & f - y_1 - \sum_{i=1}^s \sigma_i g_i \in \SOS \\
		\quad & \deg(\sigma_i g_i) \leq k,	&& \quad i = 1, \ldots, s \\
		\quad & \sigma_i \in \SOS, 		&& \quad i = 1, \ldots, s \\
		\quad & y_1 \in \R.
\end{aligned}
\tag{$\textsf{\textup{Q}}_k$}
\end{equation}
In \ref{eq:constrained_optimization_sos}, the decision variables are $y_1 \in \R$ and the real coefficients of 
$\sigma_1,\ldots,\sigma_s \in \PP$. We then have that every feasible solution $y_1$ to \ref{eq:constrained_optimization_sos} gives a lower bound on
(\ref{eq:constrained_optimization}), i.e., on $\min\{f(x) \ | \ x \in K_S\}$.
Indeed, if $y_1$ is feasible, there are 
$\sigma_0, \ldots, \sigma_s \in \SOS$, $\deg(\sigma_i g_i) \leq k$ for
$i=1,\ldots,s$, such that 
\begin{align*}
f			& = y_1 + \sigma_0 + \sum_{i=1}^s \sigma_i g_i  \label{eq:constrained_optimization_feasible_solution} \\ 
\Longrightarrow \ f(x) 	& = y_1 + \sigma_0(x)
	+ \sum_{i=1}^s \sigma_i(x) g_i(x)  \geq y_1, \quad x \in K_S, 
\end{align*}
as $\sigma_i \in \SOS$, hence $\sigma_i$ are nonnegative, and $g_i(x) \geq 0$ on $K_S$ by definition of $K_S$.
Hence $f$ is bounded from below by $y_1$ on $K_S$, i.e., every feasible solution to \ref{eq:constrained_optimization_sos} is a lower bound on (\ref{eq:constrained_optimization}).
A justification for the ansatz \ref{eq:constrained_optimization_sos} is the following well-known and easy consequence of Putinar's Positivstellensatz (Theorem \ref{thm:putinar_pss}):
\begin{cor}
\label{cor:approximating_sequence}
Let $M_S$ be Archimedean. Denote the minimum of 
(\ref{eq:constrained_optimization}) by $f^*$ and the minimum of \ref{eq:constrained_optimization_sos} by $y_1^{(k)}$. Then $y_1^{(k)} \nearrow f^*$ for $k \rightarrow \infty$.
\end{cor}

Although finite convergence is not guaranteed~\cite{lasserre2001global}, there are cases where an optimal solution $x \in K_S$ to (\ref{eq:constrained_optimization}) can be extracted from \ref{eq:constrained_optimization_sos}, see e.g.~\cite{henrion2005detecting}. In the unconstrained 
case $\min\{ f(x) \ | \ x \in \R^n\}$ given by $s=0$ (i.e. $K_S=\R^n$) in 
(\ref{eq:constrained_optimization}) even more is
known: Instead of solving \ref{eq:constrained_optimization_sos} with respect to $K_S=\R^n$ which would be given as $\max\{y_1 \ | \ f-y_1~\in~\Sigma\}$, one can consider the gradient variety\footnote{In case of unconstrained continuous minimization, provided minimizers exist, restricting minimization of $f$ to the subset of $\R^n$ where the gradient vanishes does not change the set of optimal solutions.}, resulting in $2n$ constraints corresponding to the equations $\partial_{x_1} f = \ldots = \partial_{x_n} f=0$ and solve
$\textsf{\textup{Q}}'_k$ with respect to 
\begin{equation}
S'=\{ \partial_{x_1} f, \ldots, \partial_{x_n} f, - \partial_{x_1} f, \ldots, - \partial_{x_n} f\}. 
\label{eq:gradient_variety}
\end{equation}
Then we have:

\begin{thm}[\cite{nie2006minimizing}]
\label{thm:generic_existence}
Consider the set of polynomials of degree at most $d \in \N_0$ that possess a global continuous minimizer: 
$$\mathcal F_d := \{ f \in \PP \ | \deg(f) \leq d \text{ and } \ \exists x^* \in \R^n \text{ s.t. } f(x^*) = f^* = \inf_{x \in \R^n} f(x) \}.$$
Then, for the sos-programs $\textsf{\textup{Q}}'_k$ with gradient variety constraints $S'$ from \eqref{eq:gradient_variety}, finite convergence holds for almost all polynomials $f \in \mathcal F_d$. More precisely, 
there is a $k_0~\in~\N_0$ s.t. for the optimal solutions $y_1^{(k)}$ of $\textsf{\textup{Q}}'_k$ one has $y_1^{(k)} = y_1^{(k_0)} = f^*$ for $k\geq k_0$. Moreover, a minimizer $x^*$ of \eqref{eq:constrained_optimization} can then be extracted.
\end{thm}

\section{Existence of minimizers: sufficient and necessary conditions}
\label{sec:existence_of_minimizers}

Before we search for integer minimizers of a polynomial $f \in \PP$, we review sufficient and necessary conditions to decide whether integer or continuous minimizers exist at all. For nonconstant univariate polynomials, this is equivalent to an even degree and a positive leading coefficient, which is in turn closely related to the behavior of $f(x)$ as $|x| \rightarrow \infty$. For multivariate $f$, the situation is similar once we decompose $f$ into its homogeneous components (see Section~\ref{sec:preliminaries} for the definition). A positive definite leading form is a sufficient condition for the existence of continuous minimizers whilst positive semidefiniteness is a necessary condition~\cite{marshall2003optimization,marshall2009representation}. In our next result we show that this holds for integer minimizers as well. 
Together with some observations that will be of use later on, these results are reorganized in the following proposition.

\begin{prp}
\label{prp:complete_characterization}
Let $f \in \R[X]$ with $\deg f=d>0$.
The following implications hold:
$$
\xymatrix@C=1em@R=1em@M+=0.3em{
	f_d > 0 \ar@{=>}[r] & \text{all }\sls^f(z) \text { compact} \ar@{=>}[r] \ar@{<=>}[d] & f \text{ has i.m. } \ar@{=>}[r] & \displaystyle \inf_{x \in \Z^n} f(x) > -\infty  \ar@{=>}[r] & f_d \geq 0\ar@{=>}[r] & d \text{ even} \\
& \displaystyle \liminf_{|x|\rightarrow +\infty} f(x) = + \infty \ar@{=>}[r] & f \text{ has c.m.} \ar@{=>}[r] & \displaystyle \inf_{x \in \R^n} f(x) > - \infty  \ar@{=>}[u] \\
}
$$
where i.m. abbreviates integer and c.m. continuous minimizers.
In addition, none of the implications above can be strengthened.
\end{prp}

\begin{proof}
	Let $f_d > 0$, and $c^*_j := \min_{x \in \sphere^{n-1}_p} f_j(x)$, $j=0, \ldots, d$, where $c^*_j > -\infty$ for $j~=~0, \ldots, d-1$ by compactness of the sphere and $c^*_d> 0$ since $f_d$ is positive definite.
	For $0 \neq x \in \R^n$ and $p \in [1, \infty]$ this means
	$$f(x) = \sum_{j=0}^d f_j (x) = \sum_{j=0}^d f_j \left( \frac{x}{\|x\|_p} \right) \|x\|_p^j \geq \sum_{j=0}^d c^*_j \|x\|_p^j.$$
	The expression on the right can be considered as a univariate polynomial in $\|x\|_p$ with positive leading coefficient, so $\liminf_{|x| \rightarrow + \infty} f(x) = +\infty$ follows. If the latter holds, the sublevel sets $\sls^f(z)$ must be bounded for all $z \in \R$. As $f$ is continuous, the sublevel sets are moreover closed, and compactness follows. In case the limit inferior is $s \in [-\infty, +\infty)$, pick $z \in (s, +\infty)$. There must be a sequence $x_k \in \R^n$, $\|x_k\|_p \rightarrow \infty$ as $k \rightarrow \infty$, such that $f(x_k) \leq z$ for all $k$. Put differently, $x_k \in \sls^f(z)$ for all $k$, hence the level set is unbounded, and the only equivalence in the diagram is proven. We show the two rightmost implications in the first row next, the remaining ones are straightforward. So suppose there is $x \in \R^n$ such that $f_d(x)<0$. By homogeneity, we may assume $x \in \sphere^{n-1}_\infty$. By continuity, there is a whole neighborhood $W$ of $x$ s.t. $f_d(y) < 0$ for all $y \in W$. As $W \cap \sphere^{n-1}_\infty \neq \emptyset$, there is a point $r \in W \cap \sphere^{n-1}_\infty$ with rational coordinates $r_i = \frac{z_i}{n_i}$, $z_i \in \Z$, $n_i \in \N$, $i = 1, \ldots, n$. Now for all $\lambda \in \R$, 
$$f(\lambda r)=\sum_{j=0}^d f_j(r)\lambda^j,$$
and since $f_d(r) < 0$, we have $f(\lambda r) \rightarrow - \infty$ as $\lambda \rightarrow \infty$.
Since $r_i = \frac{z_i}{n_i}$, $i=1, \ldots, n$, there is a lowest common denominator $l \in \N$ of the $r_i$. 
For $k \in \N$, we have especially $f(klr) \rightarrow - \infty$ as $k \rightarrow \infty$. But since $klr \in \Z^n$, $f$ is unbounded from below on $\Z^n$.
For the last implication, let $d$ be odd. As $f_d$ is a nonzero polynomial, there is $x \in \R^n$ s.t. $f_d(x)\neq 0$. 
Homogeneity of (odd) order $d$ implies either $f_d(x)<0$ or $f_d(-x)<0$, therefore $f_d$ is not positive semidefinite. 
\medskip

Finally, a collection of counterexamples proving that none of the implications of the proposition can be strengthened can be found in~\cite{behrends2013lower}.
\end{proof}

In the following, we rely on the sufficient condition $f_d > 0$ to ensure the existence of integer minimizers. As deciding nonnegativity of $f_d$ is NP hard, we compute a lower bound $c_d$ on the leading form $f_d$ restricted to the sphere, i.e.,  

\begin{equation}
	c_d \leq c_d^* = \min_{x \in \sphere^{n-1}_p} f_d(x).
	\label{eq:lower_bound}
\end{equation}
If $c_d >0$ we know from (\ref{eq:definiteness_by_bounds}) that $f_d >0$, so integer minimizers exist by Proposition~\ref{prp:complete_characterization}.
Our approach fails if $c_d \leq 0$ unless we find a point $x \in \sphere^{n-1}_p$ s.t. $f_d(x) < 0$ which certifies that $f_d$ is not positive semidefinite, and hence $f$ cannot have minimizers.

\section{Norm bounds on the minimizers}
\label{sec:norm_bound}

\subsection{A new bound on the norm of integer minimizers}
\label{subsec:the_bound}

If $f = \sum_{\alpha} a_\alpha X^\alpha$ satisfies $c_d^{*} = \min_{x \in \sphere^{n-1}_p} f_d(x) > 0$, i.e., $f_d$ is positive definite and a lower bound $c_d$ on the minimum with $0 < c_d \leq c_d^*$ is known, it is possible to give a bound $R \geq 0$ on the norm of the continuous minimizers. We only found one bound in the literature, which assumes $p=2$,
	\begin{equation}
	R_{\text{lit}} : = \max \left ( 1, \frac{1}{c_d} \sum_{j=1}^{d-1} \| f_j \|_1 \right ) = \max \left (1, \frac{1}{c_d} \sum_{0 < | \alpha | < d} |a_\alpha| \right ),
	\label{eq:marshall_bound}
\end{equation}

from Marshall where $\|f\|_1 := \sum_{\alpha}|a_\alpha|$ for $f = \sum_{\alpha}a_\alpha X^\alpha$;
 it is a special case (empty constraint set) of a more general result~\cite{marshall2003optimization}.  Laurent~\cite{laurent2009sums} gives a more elementary proof for Marshall's bound (\ref{eq:marshall_bound}) by showing $f(x) > f(0)$ for $\| x \|_2 > R_{\text{lit}}$. Hence $R_{\text{lit}}$ gives a valid bound on integer minimizers as well.
However, for non-sparse polynomials, this bound may get quite large. Within branch and bound approaches it is crucial to find a small bound $R$ to reduce the number of feasible solutions --  scaling $R$ by a constant $C > 0$, the number of integer points that satisfy the norm bound scales with a factor of (roughly) $C^n$. We hence 
suggest a different approach: In the following theorem, we still compute $R\geq 0$ with $f(x) > f(0)$ for $\|x\|_p>R$, but instead of bounding all homogeneous components simultaneously, we compute constants $c_j$ such that $c_j \leq c_j^* = \min_{x \in \sphere^{n-1}_p} f_j(x)$ on a suitable sphere $\sphere^{n-1}_p$. 
\begin{thm}
	\label{thm:bound_on_norm}
	Let $f \in \PP$ with $\deg f = d > 0$. For a fixed $p \in [1, \infty]$, let $c_j  \in \R$ s.t.
	$f_j(x) \geq c_j$ for all $x \in \sphere^{n-1}_p$, $j = 1, \ldots, n$. 
	Suppose $c_d > 0$. Let $R$ denote the largest nonnegative real root of the univariate polynomial $q:\R \to \R$, 
 $$q(\lambda) := \sum_{j=1}^d c_j \lambda^j.$$
	\begin{enumerate}
		\item 		Then, integer as well as continuous minimizers $x'$ of $f$ (do exist and) satisfy $\| x' \|_p \leq R.$ \label{thm:bound_on_norm:both}
	\end{enumerate}
	Let $x^*$ be any of the integer minimizers.
	\begin{enumerate}
		\setcounter{enumi}{1}
		\item We have $| x^*_i| \leq \lfloor R \rfloor$, for $i=1, \ldots, n$. \label{thm:bound_on_norm:infty}
	\end{enumerate}
	\begin{proof}
		We prove \ref{thm:bound_on_norm:both}., the other assertion follows directly from integrality of $x^*$.
		By compactness of the sphere, every $f_j$ is bounded below by some $c_j \in \R$. 
		We observed in (\ref{eq:definiteness_by_bounds}) that $c_d > 0$ implies positive definiteness of $f_d$, hence integer and continuous minimizers exist and $d$ is even (Proposition \ref{prp:complete_characterization}). Using homogeneity,
	\begin{equation}
		\label{eq:f_vs_f0}
		f(x) - f(0) = \sum_{j=1}^d f_j(x) = \sum_{j=1}^{d} f_j \left ( \frac{x}{\| x\|_p} \right ) \| x \|_p^j \geq \sum_{j=1}^d c_j \| x\|_p^j = q\left( \|x\|_p \right)
	\end{equation}
		for $x \neq 0$. Since $q$ is univariate and of degree $d > 0$, it has at most $d$ real roots. As $q(0) = 0$, $q$ has roots in $[0, \infty)$, and we denote the largest of them by $R$. As before, $c_d > 0$ yields $\lim_{\lambda \rightarrow +\infty}q(\lambda) = +\infty$. This together with the intermediate value theorem implies $q(\lambda) > 0$ for $\lambda > R$. Thus, eq. (\ref{eq:f_vs_f0}) forces $f(x) > f(0)$ for $\|x\|_p > R$.
	\end{proof}
\end{thm}

\begin{rem}
\label{rem:monoton}
The larger the $c_j$ the smaller the resulting 
norm bound $R$. Formally, let $q = \sum_{j=1}^n c_j \lambda^j$, $\tilde q = \sum_{j=1}^n \tilde c_j \lambda^j$, such that $c_j \geq \tilde c_j$, and call the largest nonnegative roots $R$ and $\tilde R$, respectively. Wlog, it suffices to consider the case that $c_j = \tilde c_j$ for $j\neq k$ and $c_k > \tilde c_k$ for some $k \in \{1, \ldots, n\}$. Now $q(\lambda) - \tilde q(\lambda) = (c_k - \tilde c_k) \lambda^k > 0$  for $\lambda > 0$ and by assumption on $c_k$, $\tilde c_k$. Thus, $q(\lambda) > \tilde q(\lambda)$ for $\lambda > 0$, hence $R < \tilde R$~--~unless $\tilde R=0$. In this case $R = \tilde R = 0$.
\end{rem}

Before we present different methods of computing valid $c_j$, we compare $R$ and $R_{\text{lit}}$. In the experiments in Section~\ref{sec:implementation}, we show that our norm bound $R$ is drastically smaller than $R_{\text{lit}}$. Also, it can be proven that our norm bounds are never larger and, except for special cases, are actually strictly smaller than the bound from the literature. 

\begin{prp}
	\label{prp:comparison_of_bounds}
	Let $f$ with $\deg f = d > 0$ and $c_d >0$ such that $f_d(x) \geq c_d$ for all $x \in \sphere^{n-1}_2$, $j = 1, \ldots, n$. 
        Compute 
	$R \in [0, \infty)$ as in Theorem \ref{thm:bound_on_norm} for 
	\begin{equation}
	c_j := - \| f_j \|_1, \quad j = 1, \ldots, d - 1,
		\label{eq:cj_brute}
	\end{equation}
	and compute $R_{\text{lit}} \in [1, \infty)$ as in (\ref{eq:marshall_bound}). 
	Then $R \leq R_{\text{lit}}$. If moreover $d >2$ and there is a coefficient $a_\alpha \neq 0$ of $f$ with $|\alpha| < d-1$, then $R < R_{\text{lit}}$ for $R \neq 1$ and $R = R_{\text{lit}}$ for $R=1$.
	\begin{proof}
		At first we observe that the numbers $c_j = - \| f_j\|_1$ in (\ref{eq:cj_brute}) are indeed valid lower bounds, for general $p \in [1,\infty]$: As $\| x\|_p \leq 1$ implies $\|x\|_\infty \leq 1$ and hence $|x^\alpha| \leq 1$, one has
$$f_j(x) = \sum_{|\alpha| = j} a_\alpha x^{\alpha} \geq \sum_{|\alpha| = j} - |a_\alpha| |x^\alpha| \geq \sum_{|\alpha| = j} -|a_\alpha| = - \| f_j\|_1 = c_j, \ x \in \sphere^{n-1}_p.$$
We prove the case $d > 2$ and $a_\alpha \neq 0$ for some $\alpha$ with $|\alpha| < d -1$. The claim obviously holds in case $R<1$. For the cases $R=1$ and $R>1$, define $q(\lambda) = \sum_{j=1}^{d} c_j \lambda^j$ as before and let $\tilde q(\lambda) = c_d \lambda^d + \left(\sum_{j=1}^{d-1} c_j \right)\lambda^{d-1}$. Then we have
\begin{equation}
	q(\lambda ) > \tilde q(\lambda) \text{ for } \lambda > 1, \quad q(\lambda) = \tilde q(\lambda) \text{ for } \lambda = 1
	\label{eq:qq-comparison}
\end{equation}
as $c_j \leq 0$ for  $j = 1, \ldots, d-1$ and one $c_k < 0$ for some $k \in \{1, \ldots, d-2\}$ by the assumption on $a_\alpha$.
By definition, the largest nonnegative real root of $q$ is $R$, and the largest nonnegative real root $\tilde R$ of $\tilde q$ is
$$\tilde R = - \frac{1}{c_d} \sum_{j=1}^{d-1} c_j = \frac{1}{c_d} \sum_{0 < | \alpha | < d} | a_\alpha|$$
and, by definition, $R_{\text{lit}} = \max(1, \tilde R)$. 
If $R = 1$, we infer from \eqref{eq:qq-comparison} that $0 = q(1) = \tilde q(1)$, so $R_{\text{lit}} = 1$.
In case $R > 1$, we infer from \eqref{eq:qq-comparison} that $0 = q(R) > \tilde q(R)$, so $R < R_{\text{lit}}$ as $\tilde q(\lambda) \rightarrow + \infty$ for $\lambda \rightarrow + \infty$.
The proof for the two remaining cases, $d = 2$ or all $a_\alpha = 0$ for $|\alpha | < d -1$, is similar as $q = \tilde q$ in these cases.
\end{proof}
\end{prp}

We now present different ways of computing bounds $c_j$ on $c_j^*=\min_{x \in \sphere^{n-1}_p} f_j(x)$.

\begin{enumerate}
	\item 
	\label{item:cj_brute:1}
	We saw in the proof of Proposition \ref{prp:comparison_of_bounds} that $c_j = - \| f_j\|_1$ gives valid lower bounds for any $p \in [1, \infty]$. However, this bound is rather rough and only useful for the lower order forms, that is those $f_j$ with $j < d$.

	\item \label{item:cj_best} The arguably easiest way to find such $c_j$ by sos programming is to minimize $f_j$ on the sphere $\sphere^{n-1}_p$: More specifically, for $p \in 2 \N$, the constraint $\| x \|_p = 1$ is equivalent to the constraint $ \sum_{i=1}^n x_i^p = 1$, which is semi-algebraic for even $p > 0$. The hierarchy \ref{eq:constrained_optimization_sos} with $g_1 = 1 - \sum_{i=1}^n X_i^p$ and $g_2 = \sum_{i=1}^n X_i^p - 1$ can be rewritten as
		\begin{equation}
		\label{eq:norm_bound_using_sos}
		\begin{aligned}
			\max 	\quad & y_1 \\
			\text{s.t.}	\quad & f_j - y_1 - q\cdot \left(1-\sum_{i=1}^n X_i^p \right) \in \Sigma \\
					\quad & q \in \PP,\ \deg q \leq k \\
					\quad & y_1 \in \R 
		\end{aligned}
		\end{equation}
		where we used $\sigma_1 g_1 + \sigma_2 g_2 = (\sigma_1 - \sigma_2)g_1 = q g_1$, some $q \in \PP$, as any polynomial can be written as the \textit{difference} of sums of squares, e.g. using $4 q = (q+1)^2 - (q-1)^2$.

			\item \label{item:cd_relaxation} A different lower bound on the leading form can be computed via the program
		\begin{equation}
			\max \ \gamma \quad \text{s.t.} \quad f_d - \gamma \cdot \sum_{i=1}^n X_i^d \in \SOS,
			\label{eq:Nie_bound}
		\end{equation}
		from \cite{nie2012sum} choosing $p = d$.
		\item We present two refined approaches of item \ref{item:cj_brute:1} in the Appendix: As a first step, we replace the underlying estimate $\|x^\alpha \|\leq 1$ by $\|x^\alpha\|\leq \|\hat x^\alpha\|$, where $\hat x$ is a maximizer of $x^\alpha$ on the sphere. In a second step, considering all orthants separately allows then to furthermore get rid of approximately half of the terms.  
\end{enumerate}

\begin{rem}
	If $p \in 2 \N$, the set $M_S$ with $S = \{1 - \sum_{i=1}^n X_i^p, \sum_{i=1}^n X_i^p - 1\}$ is Archimedean.
	Hence, from Corollary \ref{cor:approximating_sequence}, the optimal objective values of (\ref{eq:norm_bound_using_sos}) converge, for $k \rightarrow \infty$, to $c_j^* = \min_{x \in \sphere^{n-1}_p} f_j(x)$  -- which are, by Remark~\ref{rem:monoton}, the best possible bounds $c_j$. 
\end{rem}

\subsection{Application to systems of polynomial equations}
\label{subsec:app}

In this section, we consider an application to systems of polynomial equations; we test our approach on random instances of polynomials in the next section. It is a common approach to solve a system of equations $g_i(x)= 0$, $i = 1,\ldots, s$, with solutions restricted to, say, $x \in \Z^n$, $\Q^n$ or $\R^n$, by minimizing $f = g_1^2 + \ldots g_s^2$ over the integers, rationals or reals, respectively. If the minimum is $0$ at some $x$, the equations have a solution at $x$; if the minimum is nonzero, there cannot be any solution.

\subsubsection{Diophantine equations}
As an example, does the system 
\begin{align*}
- 3x_1^3 + x_1^2x_2 - x_1^2 + 2x_1x_2 + x_1 - 2x_2^2 - 2x_2 + 4 &= 0\\
2x_2^3 + x_1x_2^2 + 4x_2 - 5 &= 0
\end{align*}

possess an integer solution?
Denote the polynomials in $\Z[X_1, X_2]$ on the left hand side in the first and second equation by $g_1$ and $g_2$, respectively, and consider $f := g_1^2 + g_2^2$. The homogeneous components of $f$ are bounded from below on $\sphere^{1}_6$ by
$$(c_1, \ldots, c_6) = ( \num{-60.49}, -13.03, -41.76, -7.85, -24.45, 2.59),$$
we found the values by solving (\ref{eq:norm_bound_using_sos}) numerically.
The univariate polynomial $q(\lambda) = \sum_{j=1}^6 c_j \lambda^j$ has only two real roots: $0$ and $R \approx 9.90$.
Thus, by Theorem \ref{thm:bound_on_norm}, integer minimizers exist and must be in the box $[-9, 9]^2$. Iterating over all integer points in the box one finds $f(x_1, x_2) = 0$ at $(x_1, x_2) = (-1, 1)$. From the perspective of number theory, our method provides \textit{search bounds} on solutions of a system of Diophantine equations if the leading form of $f = \sum_{j=1}^s g_j^2$ is positive definite. 

\subsubsection{Bounds on algebraic varieties}

Similarly to the systems of Diophantine equations, our bounds apply to real algebraic varieties: Given $g_1, \ldots, g_s \in \PP$, the variety of the $g_i$ is $V(g_1, \ldots, g_s) = \{x \in \R^n \ | \ g_1(x) = 0, \ldots, g_s(x) = 0\}$. If the leading form of $f=\sum_{j=1}^s g_j^2$ is positive definite, we may give a norm bound on all points of the variety.  
As an example, let us consider the system from~\cite[Example 2, Sec. 2 § 8]{cox2007ideals}:
\begin{align}
	x^2 + y^2 + z^2 & = 1 \label{eq:example_system_2_norm}\\
x^2 + z^2 &= y \nonumber \\
x &= z \nonumber \\
x, y, z & \in \C \nonumber 
\end{align}
Computing the $c_j$ by solving (\ref{eq:norm_bound_using_sos}) for $p=2$ 
yields $(c_1, \ldots, c_4) = (0, -2.0, -0.77, 1.0)$ 
and gives us $R\approx 1.86$ as a 2-norm bound on all points in the variety. 
It is known that the variety consists of exactly four points: The system has two real and two complex solutions $(x,y,z_i)$ with $z_i \in \{\pm \frac12 \sqrt{\pm \sqrt5 -1}\}$, where the real solutions suffice $\|(x, y, z)\|_2=1$ by (\ref{eq:example_system_2_norm}). 
We conclude that in this case our bound is not far off.

\section{A class of underestimators}
\label{sec:underestimators}

\subsection{Global underestimation}
Now let $f, g: \R^n \rightarrow \R$. We then have 
\begin{equation}
	\left ( \forall x \in \R^n: g(x) \leq f(x) \right ) \Longrightarrow \inf_{x \in \Z^n} g(x) \leq \inf_{x \in \Z^n} f(x) 
	\label{eq:underestimation}
\end{equation}
where $\inf_{x \in \Z^n} g(x)$ gives a stronger bound on the integer minimum of $f$ than $\inf_{x \in \R^n}g(x)$. Using the integer minimum of $g$ to derive a lower bound on the integer minimum of $f$ makes only sense if
integer minimization of $g$ is easy compared to integer minimization of $f$.
We motivate our class of easy-to-minimize underestimators $g$ with an observation on monomials with a shift in the argument which shall serve as the building blocks to the more general underestimators.

\begin{obs}
\label{obs-shifted}
	For some $h \in \R^n$ and $\alpha \in \N_0^n$, let
	$$g = (X - h)^\alpha = \prod_{j=1}^n (X_j - h_j)^{\alpha_j}$$
        be a \emph{shifted monomial}.
	If all $\alpha_i$ are even, $g$ has a continuous minimizer at $h$ and an integer minimizer at $\rd{h}$. 
If one $\alpha_i$ is odd, $g$ is not bounded from below and does not have continuous or
integer minimizers. 
\end{obs}

Our underestimators are conic combination of shifted monomials with even $\alpha_j$, $j=1,\ldots,n$, as the combinations inherit the integer minimizer $\rd{h}$. More precisely:

\begin{prp}
	\label{prp:our_class}
Let a polynomial $g \in \PP$ be given as $g = \sum_\alpha b_\alpha (X-h)^{2\alpha}$ with $b_\alpha \geq 0$ for $\alpha \neq 0$, and $h\in \R^n$. 
	\begin{enumerate}
		\item The restriction of $g$ to $\prod_{i=1}^{k-1} \{x_i\} \times \R \times \prod_{i=k+1}^n \{x_{i}\} $ 
			that is, the univariate function $y \mapsto g(x_1, \ldots, x_{k-1}, y, x_{k+1}, \ldots, x_n)$ for fixed $x \in \R^n$ is nonincreasing for $y \leq h_k$ and nondecreasing for $y \geq h_k$, $k \in \{1, \ldots, n\}$.
		\item We have $g(x_1, \ldots, x_n) \geq g(x_1, \ldots, x_{k-1}, \rd{h_k}, x_{k+1}, \ldots, x_n)$ for every $x \in \Z^n$.
		\item $h$ is a continuous and $\rd{h}$ an integer minimizer of $g$.

	\end{enumerate}
	\begin{proof}
		The claimed properties hold for every term $(X-h)^{2\alpha}$. Thus they hold 
for conic combinations of such terms.
	\end{proof}
\end{prp}

Three properties make these polynomials $g$ useful underestimators: integer minimization is trivial, and all nonlinearity is confined to the parameter $h$. Also, the fact that the expression is linear in the $b_\alpha$ makes them accessible to optimization.
Proposition \ref{prp:our_class} motivates 
\begin{nota}
	We denote the set of conic combinations of monomials with a shift of $h$ by
	$$\cone(h) := \left \{g \in \PP \ \bigg | \ g = \sum_{\alpha \in J} b_\alpha (X-h)^{2\alpha},\ b_\alpha \in \R_{\geq 0} \mbox{ for all } \alpha\not=0, J\subset \N_0^n \text{ finite} \right \}.$$
\end{nota}

As an example, the polynomial
$$g = (X_1 - 1.5)^4(X_2 - 2)^{6} + 0.3(X_1 - 1.5)^2(X_3 - 3.2)^{8} -1 \in \cone(1.5, 2, 3.2 )$$
with $J=\{(4,6,0),(2,0,8),(0,0,0)\}$ has an integer minimizer at $(1, 2, 3)$. 

\begin{prp}
\label{prp:bound}
Let $g \in \cone(h)$ satisfy $g(x) \leq f(x)$ for all $x \in \R^n$. Then
$$ g(\rd{h}) \leq \inf_{x \in \Z^n} f(x) $$
\begin{proof}
This follows from (\ref{eq:underestimation}) and Proposition~\ref{prp:our_class}.
\end{proof}
\end{prp}

For determining an underestimator $g$ we still have to choose $h$ and the coefficients 
$b_{\alpha}$. This is described next.

	\paragraph*{Choice of $h$:} In principle, every $h \in \R^n$ may be chosen. Heuristically, we chose an (approximate) continuous minimizer of $f$ since $g$ has its continuous minimizer at $h$. In fact, every nontrivial $g$ looks like an elliptic paraboloid or a parabolic cylinder near $h$, as does $f$ near every local minimum. 
For almost all $f$, the continuous minimizer of $f$ can be found using sos methods (Theorem~\ref{thm:generic_existence}).

\paragraph*{Choice of $b_\alpha$:} 
We choose the $b_\alpha$ so that the lower bound $g(\rd{h})$ is maximized. In other words, we wish to maximize the expression
$$g(\rd{h})  = \sum_{\alpha \in J} b_\alpha (\rd{h}-h)^{2\alpha}$$
subject to $g \leq f$.
The higher order terms in $g$ ensure a certain aggressiveness in the growth behavior away from $h$, even for small coefficients $b_\alpha$, which leads to strong bounds.

Using the notation $w_\alpha := (\rd{h}-h)^{2\alpha}$, we get the following optimization problem:
\begin{align*}
\max_{J,\ b_\alpha}	\quad & \sum_{\alpha \in J} w_\alpha b_\alpha \\
\text{s.t.}		\quad & f(x)-\sum_{\alpha \in J} b_\alpha (x-h)^{2\alpha} \geq 0 \quad \forall x \in \R^n \\
			\quad & b_\alpha 			\geq 0 \quad \text{for } \alpha \neq 0 
\end{align*}
with decision variables $b_\alpha \in \R$, $\alpha \in J$ and $J \subset \N_0^n$ finite. Since this program is not tractable in general, we consider the following sos version instead:
\medskip

\begin{mdframed}
\begin{align}
y = \max	\quad & \sum_{\alpha \in J} w_\alpha b_\alpha \tag{\textup{GLOB}} \label{GLOB} \\
\text{s.t.}		\quad & f - \sum_{\alpha \in J} b_\alpha (X-h)^{2\alpha} & & \quad \text{is sos in } \R[X_1, \ldots, X_n], \nonumber \\
\quad & b_\alpha 	& & \quad \text{is sos in } \R[X_1, \ldots, X_n] \text{ for }\alpha \neq 0. \nonumber
\end{align}
\end{mdframed}
\medskip

The decision variables are the real $b_\alpha$, $\alpha \in J$. Note that $b_\alpha \in \SOS$ is equivalent to $b_\alpha \geq 0$. Once $J$ is fixed, \ref{GLOB} is a valid sos program. We show in Corollary \ref{cor:degree_bounds} that it is sufficient to choose $J = \{ \alpha \in \N_0^n \ | \ | \alpha| \leq \deg (f) /2 \}$. 
\medskip

In the following we identify a solution $b_{\alpha}$, $\alpha \in J$, with the polynomial $g$ it defines, that is with ${g = \sum_{\alpha \in J} b_\alpha \left (X-h \right)^{2\alpha}}$, and hence may say that a polynomial is a feasible or optimal solution to \ref{GLOB}.
We note that every feasible solution to \ref{GLOB} (for any choice of $h$) gives valid lower bounds on \ref{eq:problem}:

\begin{thm}
	\label{obs:properties_glob}
	Let $f \in \PP$, $h \in \R^n$ and $g = \sum_{\alpha \in J}b_\alpha (X-h)^{2\alpha} \in \cone(h)$ be a feasible solution to \ref{GLOB} for some $J$. Then 
	\begin{enumerate}
		\item $g(\rd{h}) \leq \inf_{x \in \Z^n} f(x)$. \label{obs:properties_glob:lb}
	\end{enumerate}
	If moreover $f - f(h) \in \SOS$ holds and $g$ is an optimal solution to \ref{GLOB}, then 
	\begin{enumerate}
			\setcounter{enumi}{1}
		\item $g(\rd{h}) \geq f(h)$. \label{obs:properties_glob:better}
	\end{enumerate}
	\begin{proof}
		Claim \ref{obs:properties_glob:lb} holds as $g$ being feasible to \ref{GLOB} implies $f-g \in \SOS$, hence $f-g \geq 0$, and the claim follows by Proposition~\ref{prp:bound}.
		Concerning Claim~\ref{obs:properties_glob:better}, observe that $f - f(h) \in \SOS$ implies that $h$ is a continuous minimizer of $f$ and that the constant polynomial $\tilde{g}=f(h)$ is a feasible solution to \ref{GLOB},
hence $ g(\rd{h}) \geq \tilde{g}(\rd{h})=f(h)$ 
for every optimal solution $g \in \cone(h)$.
	\end{proof}
\end{thm}

\subsection{Improving the underestimators}

\subsubsection*{Motivation}

A quite restrictive condition in \ref{GLOB} is that it requires $g(x) \leq f(x)$ globally, i.e., for all $x \in \R^n$. Actually, this is not necessary for our purposes. It is enough to require $g(x) \leq f(x)$ only for those $x \in \R^n$ that satisfy $f(x) \leq f(q)$ for some $q \in \Z^n$. That is, for all $q \in \Z^n$, we have
\begin{equation}
 \left( \forall x \in \sls^{f}(f(q)): g(x) \leq f(x) \right ) \Longrightarrow \inf_{x \in \Z^n} g(x) \leq \inf_{x \in \Z^n} f(x),
	\label{eq:underestimation_sublevelset}
\end{equation}
in other words, the integer minimum of $g$ is a lower bound on the integer minimum of $f$ even if $g$ is an underestimator of $f$ only on a sublevel set $\sls^{f}(f(q))$. If we make use of this in our sos program, the lower bound can only improve.

\begin{figure}[h]
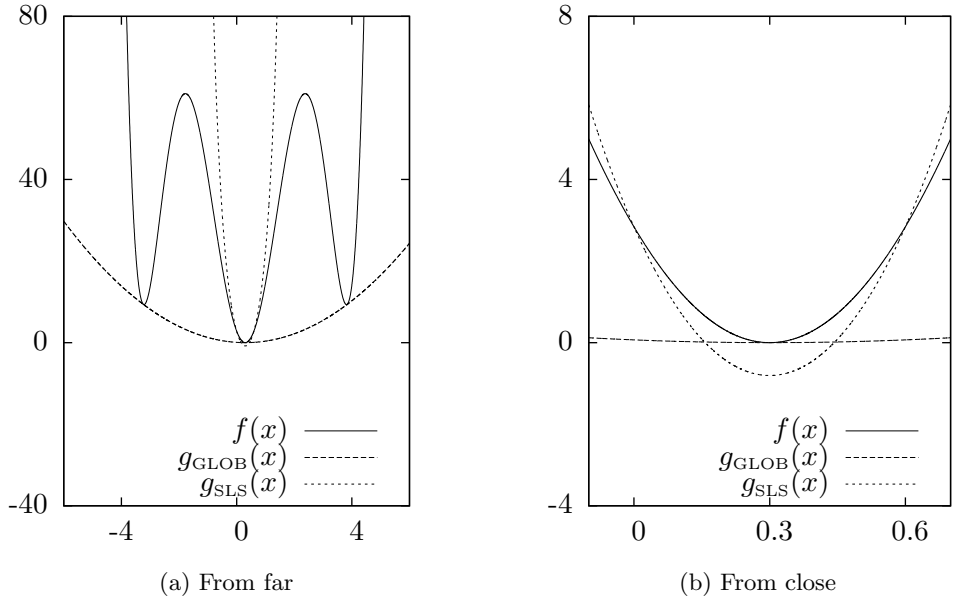

	\centering
	\subfloat[From far]{\input{plotting/uest_illustration/underestimation_from_far.tex} \label{fig:big_picture} } \quad 
	\subfloat[From close]{ \input{plotting/uest_illustration/underestimation_from_close.tex} \label{fig:closer_look}}
	\caption{Global underestimator $g_{\text{GLOB}}$ and an underestimator $g_{\text{SLS}}$ on a sublevel set.}
	\label{fig:uest_illustration}
\end{figure}

\medskip
But before we delve into the details, let us consider the potential payoff by taking a look at the example in Figure \ref{fig:big_picture}.  The plot depicts the univariate polynomial
$$f = 0.2 \cdot (X - 0.3)^6 - 5\cdot(X - 0.3)^4 + 32\cdot (X - 0.3)^2.$$
along with two underestimators $g_{\text{GLOB}}$, $g_{\text{SLS}}$.
A short calculation shows that $f$ has five local extrema at $0.3$ and $0.3 \pm \sqrt{\frac{25 \pm \sqrt{145}}{3}}$, and that 
the local minimizers are at $x=0.3$ and at $x_{\pm} = 0.3 \pm \sqrt{\frac{25 + \sqrt{145}}{3}} \approx 0.3 \pm 3.51$. Considering that $f$ has a positive definite leading form, one of the local minimizers must be a global one, and comparing the function values shows that $x=0.3$ is the continuous minimum. Moreover, $f$ must have its integer minimizer in $[-3,3]$ as $\min \{f(x_+), f(x_{-})\} > f(0) $; comparing the function values shows that $f$ has a single integer minimizer at $x = 0$ with value $f(0) \approx 2.84$. The underestimator $g_{\text{GLOB}} \in \cone(h)$, computed as optimal solution to \ref{GLOB} is given by\footnote{For this example we solved \ref{GLOB} for $h = 0.3$ and $\deg g = 6$, using SOSTOOLS 3.00 and CSDP 6.1.0.}
$$g_{\text{GLOB}} \approx   8.71 \cdot 10^{-11} \cdot (X-0.3)^6 + 1.09 \cdot 10^{-09} \cdot (X-0.3)^4 + 0.75 \cdot (X-0.3)^2 - 1.22 \cdot 10^{-09},$$
is globally below $f$. To find an underestimator on a sublevel set, we first fix the level $z=f(q)$ heuristically.
Note that any $q \in \Z$ is a feasible solution to \ref{eq:problem} and hence an upper bound; any integer minimum must be contained in $\sls^f(f(q))$. 
 As $h=0.3$ is the global minimizer, we choose $q = \rd{h} = 0$ here. The polynomial $g_{\text{SLS}}$, given by  
$$ g_{\text{SLS}} \approx   9.09 \cdot (X-0.3)^6 + 11.80 \cdot (X-0.3)^4 + 39.36 \cdot (X-0.3)^2 - 0.81,$$
is an underestimator on the sublevel set $\sls^{f}(f(0)) = [0, 0.6]$, as can be seen in Figure \ref{fig:closer_look}. It will be shown in the next section how this function can be found.
The plot reveals the shortcomings of global underestimation: Any \textit{global} underestimator in $\cone(0.3)$ cannot go above the local minimizers of $f$. This ``barrier'' from above turns $g_{\text{GLOB}}$ in this example essentially into a quadratic underestimator for small $x$ as the ratio of the higher order coefficients and the one in front of the quadratic term is of order $10^{-10}$. 
The underestimator $g_{\text{SLS}}$ however is a degree 6 polynomial whose higher order coefficients are not small at all.
Note that $g_{\text{GLOB}}$ is much closer to $f$ near $0.3$ compared to the new underestimator $g_{\text{SLS}}$.
However, the quality of the resulting lower bound depends on the function values at $0$ and there 
$g_{\text{SLS}}$ is closer to $f$ than $g_{\text{GLOB}}$.
The lower bounds the two underestimators provide are $g_{\text{GLOB}}(0) \approx 0.07$ 
and $g_{\text{SLS}}(0) \approx 2.84$. In this case, we are lucky as the lower bound on the integer minimum 
and $f(0)$ coincide, showing once more that $f$ has its integer minimizer at $0$. 

\medskip
\subsubsection*{The sos program for computing the improved underestimator}
How do we compute the improved underestimator? At first, we observe that every sublevel set $\sls^f(z)$, $z \in \R$, of $f$ is semi-algebraic. Indeed, with the notation from (\ref{eq:KS}) and $\tilde S:=\{z-f\}$, we have 
$$\sls^f(z) = \{ x \in \R^n \ | \ z - f(x) \geq 0\} = K_{\tilde S}.$$
Moreover, $\sls^f(z)$ is compact if the leading form of $f$ is positive definite (see Proposition \ref{prp:complete_characterization}). Compactness of $\sls^{f}(z)$ in turn implies that the quadratic module $M_S \subset \R[X_1, \ldots, X_n]$ generated by $S:= \{f-g, z - f\}$, for any $g \in \R[X_1, \ldots, X_n]$, is thus by definition Archimedean. Hence, for every feasible underestimator $g \in \cone(h)$ the existence of a representation for $f-g$ as in Putinar's Positivstellensatz (Theorem \ref{thm:putinar_pss}) is guaranteed. This motivates the following program:
\medskip

\begin{mdframed}
\begin{align}
	y^{(k)} = \max \quad & \sum_{\alpha \in J} w_\alpha b_\alpha \label{SLS} \tag{\textup{SLS}} \\
\text{s.t.}		\quad & f - \sum_{\alpha \in J} b_\alpha (X-h)^{2\alpha} - \sigma(z - f) & & \quad \text{is sos in } \R[X_1, \ldots, X_n], \nonumber \\
\quad & b_\alpha \text{ for } \alpha \neq 0,\ \sigma 	& & \quad \text{are sos in } \R[X_1, \ldots, X_n], \nonumber \\
\quad &	\deg \sigma \leq k \nonumber
\end{align}
\end{mdframed}
\medskip

The decision variables are the real $b_\alpha$ as for \ref{GLOB} and, additionally, the real coefficients of the polynomial $\sigma$. As before, we use the notation $w_\alpha := (\rd{h}-h)^{2\alpha}$. 
\ref{SLS} is a valid sos program once $J$ and the degree of $\sigma$ are 
fixed.

\begin{thm}
	\label{thm:properties_sls}
	Let $f \in \PP$, $h \in \R^n$ and $g \in \cone(h)$ be a feasible solution to \ref{SLS} with $z \geq f(q)$ for some $q \in \Z^n$. 
	\begin{enumerate}
		\item Then $g(\rd{h}) \leq \inf_{x \in \Z^n} f(x)$. \label{thm:properties_sls:lb}
		\item If $J$ is fixed, $y^{(-\infty)} \leq y^{(0)} \leq y^{(2)} \leq y^{(4)} \leq \ldots$\footnote{Note that every sos polynomial $\sigma \neq 0$ has even degree.}  \label{thm:properties_sls:inc}
		\item If $f_d$ is positive definite, there is $k_0 \in \N_0$ such that SLS is feasible for all $k \geq k_0$. \label{thm:properties_sls:pd}
		\item \ref{SLS} with $k = -\infty$ is \ref{GLOB}. \label{thm:properties_sls:equiv} 
		\item If $f - f(h) \in \SOS$ and $g$ is optimal, then $g(\rd{h}) \geq f(h)$. \label{thm:properties_sls:better}
	\end{enumerate}

	\begin{proof} 
		Statement \ref{thm:properties_sls:lb} holds as $g$ feasible implies $f-g - \sigma(z-f) \in \SOS$. Hence $f(x)-g(x) \geq 0$ for those $x$ with $f(x) \leq z$, especially for those $x$ with $f(x) \leq f(q)$ as $f(q) \leq z$ by assumption. The claim follows by \eqref{eq:underestimation_sublevelset}.

		Statement~\ref{thm:properties_sls:inc} is clear as we only allow more coefficients for $\sigma$. 

To see Statement~\ref{thm:properties_sls:pd}, note that $\sls^f(z)$ is nonempty as $z\geq f(q)$ and moreover compact (Theorem \ref{prp:complete_characterization}), so $f(x) > c$ for some $c \in \R$ and all $x \in \sls^f(z)$. Hence $f-c \in M_{\{z-f\}}$ by Putinar's Positivstellensatz (Theorem \ref{thm:putinar_pss}). This means $f-c = \sigma_0 + \sigma(z-f)$ for some sos $\sigma_0,\sigma \in \PP$. Thus $g:= c$ is a feasible solution, and $k_0 := \deg \sigma$.

To see Statement \ref{thm:properties_sls:equiv}, we note that $k= -\infty$ corresponds to $\sigma = 0$, in which case \ref{SLS} is \ref{GLOB}.

		Statement \ref{thm:properties_sls:better} is a consequence of Statements~\ref{thm:properties_sls:inc} and \ref{thm:properties_sls:equiv} and Theorem~\ref{obs:properties_glob}. 
	\end{proof}
\end{thm}

We have not yet addressed the degree of $g$ in \ref{GLOB} and \ref{SLS} nor the degree of $\sigma$ in \ref{SLS}. The following proposition shows that once the degree of $\sigma$ in \ref{SLS} is fixed, the degree of $g$ in any feasible solution is bounded from above in terms of $\deg f$ and $\deg \sigma$. 

\begin{prp}
\label{prp:degree_bounds}
Let $f \in \PP$, $g \in \cone(h)$ with $\deg f > 0$, $\deg g > 0$, $z \in \R$ and $\sigma \in \SOS$ such that
\begin{equation}
f-g - \sigma(z-f) \text { is sos}.
\label{eq:degreeBounds}
\end{equation}

Then
$$\deg (g) \leq \deg (f) + \max\{ \deg(\sigma), 0 \}.$$

\begin{proof}
Eq. (\ref{eq:degreeBounds}) is equivalent to $f-g - \sigma (z-f) = \sigma_0$ for some $\sigma_0 \in \SOS$, or 
\begin{equation}
g+\sigma_0=f(1+\sigma) - z \sigma.
\label{eq:degreeBounds2}
\end{equation}
\begin{alignat*}{1}
\text{Hence}\quad \deg (g)& \leq \max \left\{ \deg (g), \deg (\sigma_0) \right\} \stacktext{(I)} = \deg(g+\sigma_0) \stacktext{(II)} = \deg \left(f(1+\sigma) -z \sigma \right) \\
	& \stacktext{(III)} = \max \left\{ \deg \left( f(1+\sigma) \right), \deg (z \sigma) \right\} \stacktext{(IV)} = \deg \left( f(1+\sigma) \right) \\
	&\stacktext{(V)} = \deg(f)+\deg(1+\sigma) \stacktext{(VI)} = \deg (f) + \max \{ \deg (\sigma), 0 \}.
\end{alignat*}
As $g - g(h) \in \SOS$ and $\deg g > 0$, equality (I) follows from Lemma \ref{lem:sos_nonvanishing}.
Equality in (II) follows from eq. (\ref{eq:degreeBounds2}). Using $\deg f > 0$, the equalities in (III) and (IV) follow from a typical degree argument: If $u, v \in \PP$, $\deg u \neq \deg v$, we have $u +v \neq 0$ and $\deg (u+v) = \max(\deg u, \deg v)$.
Equality in (V) holds as the degree is multiplicative, (VI) follows easily if one distinguishes the cases $\sigma = 0$, $\sigma \in \R_{\geq 0}$ and $\deg \sigma > 0$.
\end{proof}
\end{prp}

\begin{cor}
	\label{cor:degree_bounds}
Let $g \in \cone(h)$ be a feasible solution for \ref{GLOB}. Then $\deg g \leq \deg f$.
\begin{proof}
Use Proposition \ref{prp:degree_bounds} with $\sigma = 0$ and the result follows from
Statement \ref{thm:properties_sls:equiv} of Theorem~\ref{thm:properties_sls}.

\end{proof}
\end{cor}

\section{Implementation and results on random instances}

\label{sec:implementation}

\subsection{Experimental setup}

To evaluate our results, we ran computer experiments: For a fixed number of variables~$n$ and an even degree~$d$, we created instances of random polynomials

\begin{equation}
f = \sum_{|\alpha| \leq d} a_\alpha X^\alpha = \sum_{|\alpha|\leq d} a_\alpha X_1^{\alpha_1} \cdots X_n^{\alpha_n}, \quad a_\alpha \sim \uniDist(-1,1) \text{ i.i.d..}
\label{eq:random_poly}
\end{equation}

\begin{figure}
	\centering
		\input{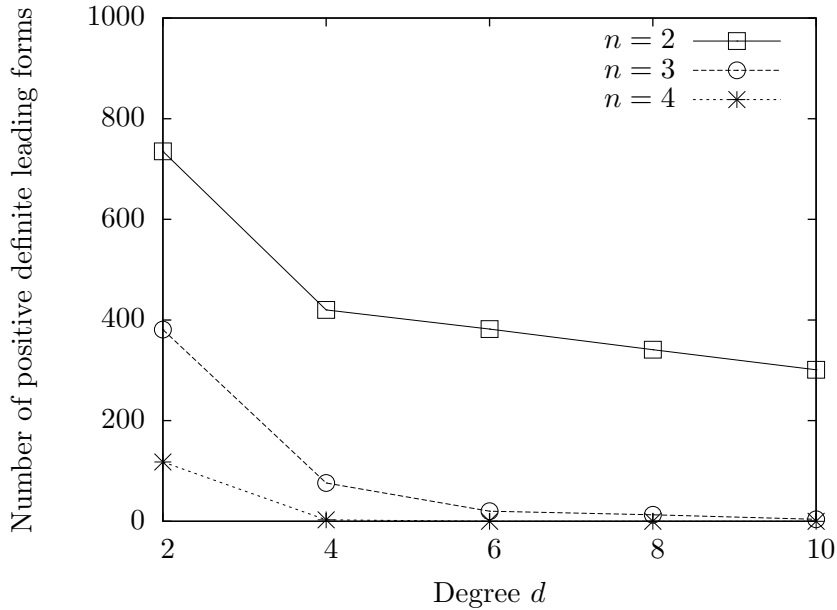}
		\caption{Instances with detected positive definite leading form.}
		\label{fig:share_of_pd_lf}
\end{figure}

As we are only interested in polynomials with positive definite leading form, we restricted ourselves 
to those polynomials that satisfy

\begin{equation}
	a_{(d,0,\ldots,0)} > 0,\ a_{(0,d,0,\ldots,0)} > 0,\ \ldots,\ a_{(0,\ldots,0,d)} > 0 \tag{A}
	\label{eq:discarding_condition}
\end{equation}

since a polynomial with at least one of these coefficients nonpositive cannot be positive definite. 
Then, we solved program~(\ref{eq:norm_bound_using_sos}) with $k =  d  + 2$ to compute a lower bound $c_d$ on $\min_{x \in \sphere^{n-1}_2} f_d(x)$ to determine whether $f$ indeed has a positive definite leading form. If $c_d \leq 0$, we discarded the instance, else we know that $f_d$ is positive definite. In the first part of the experiments, for every tuple $(n,d)$ with $n=2, 3, 4$ and $d=2, 4, 6, 8, 10$, we created 1000 random instances of polynomials that satisfy condition (\ref{eq:discarding_condition}). In Figure~$\ref{fig:share_of_pd_lf}$ we plot how many of these have been detected to satisfy $f_d > 0$. As $d$ and $n$ increase, the probability of positive definiteness should decrease -- as, loosely speaking, more (independent) random variables $a_\alpha$ simultaneously influence the result -- which is reflected in the plot. We then use these instances to evaluate the norm bounds (see Section~\ref{subsec:evaluting_norm_bounds}). In the second part of the experiments, for four tuples $(n,d)$, we again generated polynomials according to \eqref{eq:random_poly} and took the first 50 of them that were detected to have a positive definite leading form as input for the optimization problem which is in turn solved by branch and bound (see Section~\ref{subsec:evaluating_underestimators}).

\medskip

We use MATLAB\footnote{MATLAB is a registered trademark of The MathWorks Inc., Natick, Massachusetts}~2014b 64-bit, SOSTOOLS~3.00~\cite{sostools} to translate the sos programs into semidefinite programs and CSDP~6.1.0~\cite{borchers1999csdp} to solve the latter. The experiments were conducted on GNU/Linux (Ubuntu 12.04) running on 2 Intel\textsuperscript\textregistered~Xeon\textsuperscript\textregistered X5650 CPUs (each 6 cores) with a total of $\SI{96}{\giga\byte}$ RAM.

\subsection{Evaluating the norm bounds}
\label{subsec:evaluting_norm_bounds}
Once positive definiteness is certificated by some $c_d > 0$, the bounds on the norm of the minimizers can be computed. We summarize the steps to compute a norm bound on the minimizer\footnote{In the algorithm, we abbreviate integer minimizer(s) to i.m. and continuous minimizer(s) to c.m..} in algorithmic form (Algorithm \ref{algo:find_the_box}).

\begin{algorithm}
	\caption{Norm bound on minimizers}

	\begin{algorithmic}
		\Input{$f \in \R[X_1, \ldots, X_n]$ with $\deg f \in 2\N$, parameters $p \in 2 \N$, $k_{\max} \in \N_{0}$}
		\State $k \gets 0$
		\State $c_d \gets -\infty$
		\State $x\gets \NULL$
		\While{$k \leq k_{\max}$ and $c_d < 0$ and $x = \NULL$}
			\State solve program (\ref{eq:norm_bound_using_sos}) for $j=d$ and parameter $k$
			\State $c_d \gets$ optimal value
			\If{optimal solution can be extracted} 
				\State $x \gets$ optimal solution
			\EndIf
			\State $k \gets k + 1$
		\EndWhile
		\If{$c_d < 0$ and $x \neq \NULL$} 
		\Output{$x$}
		\Print{$f_d(x) < 0 $ so $f$ has neither i.m. nor c.m..} \Comment{Proposition \ref{prp:complete_characterization}}
		\ElsIf{$c_d \leq 0$}
		\Print{Cannot decide $f_d > 0$ for $k\leq k_{\max}$.}
			\Else \Comment{$c_d > 0$ in the following}
			\Print{$f$ has integer and continuous minimizers.} \Comment{$f_d > 0$ by (\ref{eq:definiteness_by_bounds})}
			\For{$j=1, \ldots, d-1$}
			\State $c_j \gets$ $\max$ of (\ref{eq:cj_brute}) and (\ref{eq:norm_bound_using_sos}) \Comment{can be improved by also taking (\ref{eq:cj_lagrange}) into account}
			\EndFor
			\State define $q: \R \rightarrow \R$, $q(\lambda) = \sum_{j=1}^d c_j \lambda^j$
			\State $R \gets$ largest root of $q$ in $\R$ \Comment{$R \geq 0$ by Theorem \ref{thm:bound_on_norm}}
			\Output{$R$}
			\Print{The minimizers $x'$ suffice $\|x'\|_p \leq R$.} \Comment{Theorem \ref{thm:bound_on_norm}}
		\EndIf

	\end{algorithmic}
	\label{algo:find_the_box}
\end{algorithm}

\begin{figure}[h]
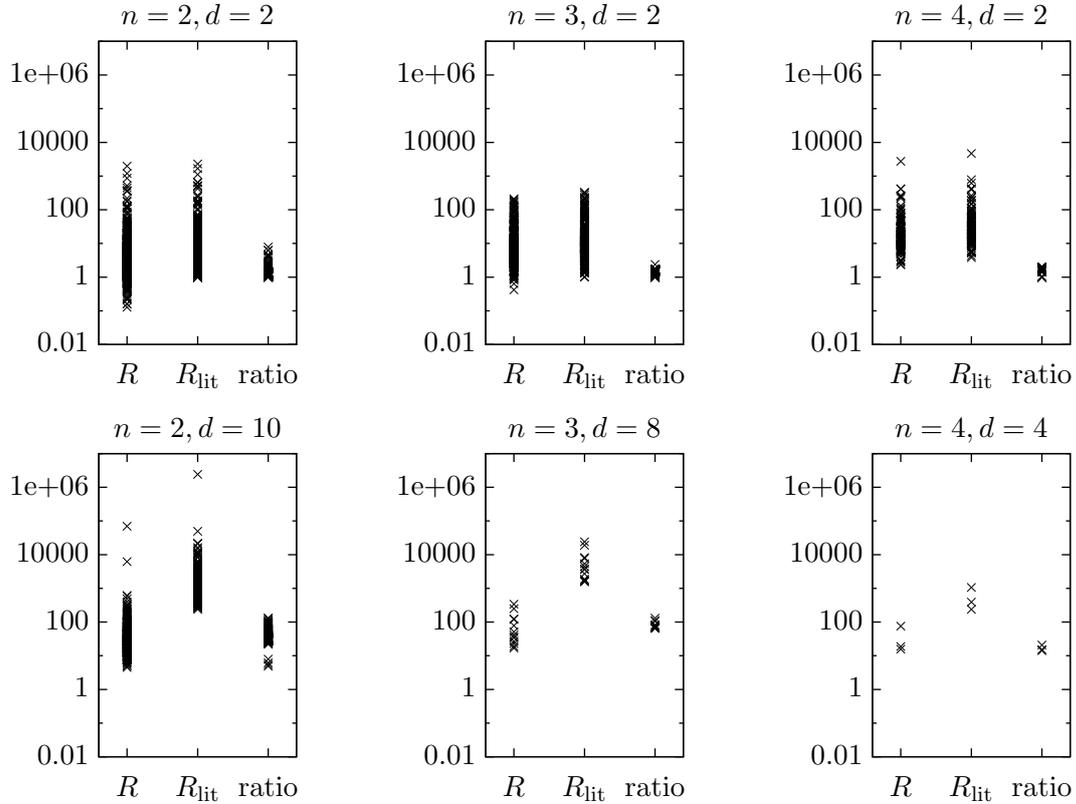

	\centering
	\subfloat{ \input{plotting/radius_experiment/radius_n_2-d_2.tex} } \qquad
	\subfloat{ \input{plotting/radius_experiment/radius_n_3-d_2.tex} } \qquad
	\subfloat{ \input{plotting/radius_experiment/radius_n_4-d_2.tex} }\\
	\subfloat{ \input{plotting/radius_experiment/radius_n_2-d_10.tex} } \qquad
	\subfloat{ \input{plotting/radius_experiment/radius_n_3-d_8.tex} }\qquad
	\subfloat{ \input{plotting/radius_experiment/radius_n_4-d_4.tex} } \vspace{0.5cm}
	\caption{Bounds on the norm of minimizers for different dimensions $n$ and degrees $d$.}
	\label{fig:norm_bounds}
\end{figure}

For $n$ and $d$ as described above, and each of the 1000 randomly created polynomials that has been detected to have a positive definite leading form, we computed the bound from the literature $R_{\text{lit}}$ from eq. (\ref{eq:marshall_bound}) and our new 
bound $R$ (Theorem \ref{thm:bound_on_norm}) such that $\|x'\|_2~\leq~\min(R_{\text{lit}}, R)$ holds for every continuous and integer minimizer of $f$.
Figure~\ref{fig:norm_bounds} depicts a selection of $d$ and $n$:
Those with smallest degree, $d=2$, and the maximal degree $d\leq10$ such that we still detected some instances with $f_d > 0$. We conclude that for the quadratic case $d=2$, our approach does not yield significantly better results. However, for a higher number of variables and $d \geq 4$ we outperform the classic norm bound on all instances. Most prominently of this selection, for $(n,d)=(3,8)$, we are better by a factor $C = R_{\text{lit}}/R$ of 50 throughout and in some instances we are better by a factor of $C \approx 100$. This means the number of feasible solutions decreases by a factor of up to $C^n \approx 100^3$ in this example.

\subsection{Evaluating our underestimators within branch and bound}

We evaluated the underestimators in a branch and bound framework. Firstly, we present an algorithm that shows how special properties of our underestimators can be exploited to speed up branching and pruning. In the actual experiments, we generated polynomials according to \eqref{eq:random_poly}, where we restricted ourselves to the tuples $(n,d) = (2, 4)$, $(2,6)$, $(3, 4)$ and $(4,2)$ to keep the problem size tractable and to have an acceptably high ratio of positive definite polynomials (compare Figure \ref{fig:share_of_pd_lf}). We generated random polynomials until we had 50 that were detected to have a positive definite leading form and which were then 
used as input to the optimization problem. In the following we present an evluation of the initial lower bound $g(\rd{h})$ and a runtime comparison with other lower bounds from the literature.

\label{subsec:evaluating_underestimators}

\subsubsection{Algorithm}

\label{sec:bandb}

Our branch and bound framework is depth first. This keeps memory usage small and allows us to quickly obtain good feasible solutions. We do not reorder the variables. Subproblems are collected in a list $\subproblems$; every subproblem $\mathcal{P} \in \subproblems$ is of the form $\mathcal{P} = (m, r_1, \ldots, r_m)$, where $m \in \{0, \ldots, n \}$ encodes the number of fixed variables $(r_1, \ldots, r_m) \in \Z^m$; i.e.,
\begin{equation}
	\begin{aligned}
		 \min 	& \quad f(r_1, \ldots, r_m, x_{m+1}, \ldots, x_n) \\
			& \quad x_{m+1}, \ldots, x_{n} \in \Z
	\end{aligned}
	\tag{$\mathcal{P} = (m, r_1, \ldots, r_m)$} 
	\label{eq:P}
\end{equation}
 and $(0)$ encodes the initial problem. Algorithm~2 states the whole procedure. 

\begin{algorithm}
	\begin{algorithmic}[4]
		\Input{$f \in \R[X_1, \ldots, X_n]$, $h \in \R^n$, $p$-norm bound $R$ on minimizers, $k \in 2 \N_0$}
		\State $x^* \gets \rd{h}$ \Comment{initial guess for integer minimizer}
		\State $u \gets f(x^*)$ \Comment{upper bound on integer minimum} 
		\State $\subproblems \gets \{ (0) \}$ \Comment{initial list of subproblems}
		\State find underestimator $g$: solve \ref{SLS} with $h$, $\deg g \leq \deg \sigma = k$ \Comment{or \ref{GLOB}, resp.}
		\While{ $\subproblems \neq \emptyset$ } 
			\State pick $\mathcal{P} = (m, r_1, \ldots, r_m) \in \subproblems $ with $m$ maximal
			\State $\subproblems \gets \subproblems \setminus \{\mathcal{P}\}$
			\If{ $m < n$ }
				\State $L \gets \left \lfloor \sqrt[\leftroot{-2}\uproot{2}p]{R^p - |r_1|^p - \dots - |r_m|^{p}} \right \rfloor$
				\State let $\tilde g : \R \rightarrow \R$, $\tilde g(x_{m+1}) =g(r_1, \ldots, r_m, x_{m+1}, \rd{h_{m+2}}, \ldots, \rd{h_{n}})$
				\If{$\tilde g(\rd{h_{m+1}}) \leq u$}  \Comment{otherwise prune}\label{code:prune}
				\State find $L_1 \in [-L, L] \cap \Z $ minimal with $\tilde g(L_1) \leq u$ \label{code:L1}
				\If{such an $L_1$ exists}
				\State find $L_2 \in [-L, L] \cap \Z $ maximal with $\tilde g(L_2) \leq u$ \label{code:L2}
				\Else \State $L_1 \gets + \infty$, $L_2 \gets  -\infty$. \EndIf
				\ForAll{$r_{m+1} \in [L_1, L_2] \cap \Z$} \label{code:branch} \Comment{$[L_1, L_2] = \emptyset$ if $L_1 = +\infty$}
					\State $\subproblems \gets \subproblems \cup \{(m+1, r_1, \ldots, r_{m+1})\}$ \Comment{actual branching}

					\EndFor
				\EndIf
			\Else \Comment{all variables $x_i$ were fixed to values $r_i$}
				\If{ $f(r) < u$ } \Comment{update upper bound}
					\State $x^* \gets r$
					\State $u \gets f(r)$
				\EndIf

			\EndIf
		\EndWhile
		\Output{$x^*$, $u$}
		\Print{$f$ attains its integer minimum $u$ at $x^*$.}
	\end{algorithmic}
	\caption{Branch and Bound}
	\label{algo:minimize}
\end{algorithm}

\begin{prp}
	\begin{enumerate}
		\item Algorithm \ref{algo:minimize} is correct, that is, it always terminates after a finite number of steps with an optimal integer solution $x^*$ that satisfies $f(x^*) = u$. \label{prp:algo_is_correct:termination}
		\item The integers $L_1$ and $L_2$ (in lines \ref{code:L1} \& \ref{code:L2}) can be found with binary search in $\lceil \log_2(L) \rceil + 2 \leq \lceil \log_2(R)\rceil + 2$ evaluations of $\tilde g$ if $L > 0$.\label{prp:algo_is_correct:steps}
	\end{enumerate}

	\begin{proof}
		Let $x^*$ be any optimal solution. To prove \ref{prp:algo_is_correct:termination}. it suffices to show that the algorithm terminates \label{lem:algo_terminates:termination} and no problem with $(n, x^*)$ as subproblem gets pruned in Step \ref{code:prune} or lost in Step \ref{code:branch}.
		To see termination of the algorithm, we observe that the number of subproblems is finite as the sets $B_m = \{ y \in \Z^m \ | \ \| y \|_p \leq R \}$, $m=1, \ldots, n$ are finite, every subproblem $(m, r_1, \ldots, r_m)$ suffices $(r_1, \ldots, r_m) \in B_m$ and no subproblem is inserted into the list $\mathcal L$ more than once.
		To see that $x^*$ does not get discarded in Step \ref{code:prune}, define
		\begin{equation}
			\tilde g(x_{m+1}) := g(x^*_1, \ldots, x^*_m, x_{m+1}, \rd{h_{m+2}}, \ldots, \rd{h_n})
			\label{eq:tilde_g}
		\end{equation}
		and suppose $\tilde g(\rd{h_{m+1}}) > u$. Hence 
		$$\tilde g(\rd{h_{m+1}}) > u \geq f(x^*) \geq g(x^*) \geq g (x^*_1, \ldots, x^*_m, \rd{h_{m+1}}, \ldots, \rd{h_n}) = \tilde g(\rd{h_{m+1}}),$$
		a contradiction, where we used the monotonicity property of $g$ (Proposition \ref{prp:our_class}) and that $g(x) \leq f(x)$ for $x \in \sls^{f}(f(q))$, a fortiori for $x \in \sls^f(f(x^*))$.
		Suppose that $x^*$ gets lost in Step \ref{code:branch}.
		Necessarily, $x^*_{m+1} < L_1$ or $x^*_{m+1} > L_2$. We derive a contradiction for $x^*_{m+1} < L_1$, the other case is identical. Observe that $x_{m+1}^* \in [-L, L]$ as every optimal integer solution satisfies $\sum_{j=1}^n |x^*_j|^p \leq R^p$, so we must have $|x_{m+1}^*| = \sqrt[\leftroot{-2}\uproot{2}p]{|x^*_{m+1}|^p} \leq \sqrt[\leftroot{-2}\uproot{2}p]{R^p - |x^*_1|^p - \ldots - |x^*_m|^p}$. As $x^*_{m+1}$ is integer, we may round down -- in other words, $x^*_{m+1} \in [-L, L]$. By definition of $L_1$ and Proposition \ref{prp:our_class}, we have $\tilde g(x^*_{m+1}) > u$ with $\tilde g$ from (\ref{eq:tilde_g}), thus, using Proposition \ref{prp:our_class} again,
		$$\tilde g(x^*_{m+1}) > u \geq f(x^*) \geq g(x^*) \geq \tilde g(x^*_{m+1}),$$
		a contradiction. \\
		We finally show that Claim \ref{prp:algo_is_correct:steps} holds. We prove the claim for $h_{k+1} \geq 0$, the proof for $h_{k+1} \leq 0$ is similar. In case $h_{k+1} > L$, $L_1$ exists if and only if $\tilde g(L) \leq u$ as $\tilde g(x_{k+1})$ is non-increasing for $x_{k+1} \leq h_{k+1}$ (by Proposition \ref{prp:our_class}); necessarily, $L_2 := L$. Using binary search on $[-L, L]$, $L_1$ can be found using at most $\lceil \log_2 (2 L) \rceil = \lceil\log_2(L) \rceil + 1$ further evaluations of $\tilde g$. In case $0 \leq h_{k+1}\leq L$, $L_1$ exists as $\tilde g(\rd{h_{k+1}}) \leq u$ in Step \ref{code:prune}. Again using binary search, $L_1 \in [-L, \rd{h_{k+1}}]$ can be found in no more than $\lceil \log_2(2 L)\rceil$ evaluations. As $\tilde g(x_{k+1}) = \tilde g(h_{k+1} - x_{k+1})$, it only needs at most one more evaluation of $\tilde g$ to find $L_2$, so we find both numbers in no more than $\lceil \log_2(L) \rceil + 2$ evaluations of $\tilde g$.
	\end{proof}
\end{prp}

\begin{rem}
	\label{rem:input_to_glob_and_sls}
	Concerning our implemenation, we chose $\deg g = \deg f$ for \ref{GLOB} and \ref{SLS} and $\deg \sigma = 2$ for \ref{SLS}. For the parameter $h \in \R^n$ we chose an (approximate) continuous minimizer computed via the SOSTOOLS function  \verb+findbound.m+ -- however, the algorithm accepts arbitrary $h \in \R^n$. We determined $R$ using Algorithm \ref{algo:find_the_box}.
\end{rem}

\subsubsection{The initial lower bound on the minimum}
\begin{figure}[!h]
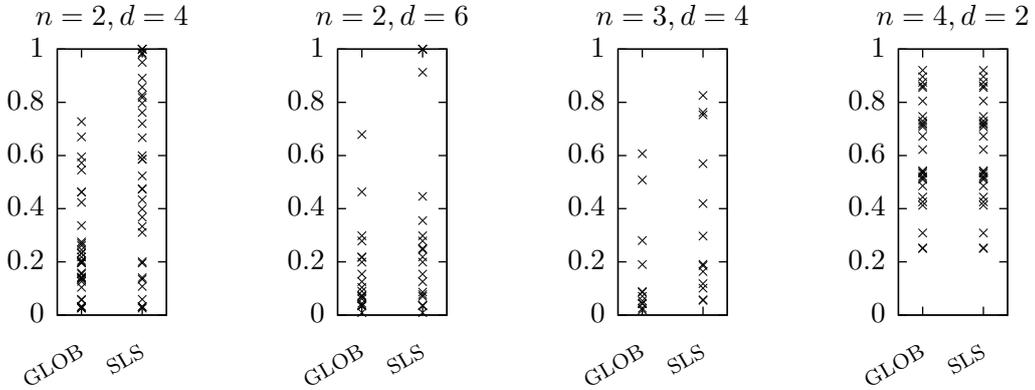

	\centering
	\subfloat{ \input{plotting/comparison_of_glob_and_sls/lower_bounds-n_2-d_4.tex} } \quad
	\subfloat{ \input{plotting/comparison_of_glob_and_sls/lower_bounds-n_2-d_6.tex} } \quad
	\subfloat{ \input{plotting/comparison_of_glob_and_sls/lower_bounds-n_3-d_4.tex} } \quad
	\subfloat{ \input{plotting/comparison_of_glob_and_sls/lower_bounds-n_4-d_2.tex} } \vspace{0.5cm} \\
	\caption{Lower bound comparison using the ratio $Q$.}
	\label{fig:lower_bound_comparison}
\end{figure}

Before we compare our underestimators with lower bounds from the literature, we directly evaluate our initial lower bound $g(\rd{h})$. To this end, we define a ratio $Q$ as follows: Let $h$ be a continuous minimizer of $f$ (if found by sos methods), $x^*$ an integer minimizer of $f$ found during B\&B and $g$ be a solution to \ref{GLOB} or \ref{SLS}.
Then
$$Q:=\frac{g(\rd{h}) - f(h)}{f(x^*) - f(h)}$$
takes values in $[0, 1]$, is invariant under scaling of $f$ by constants $\lambda > 0$ and addition of constants $c \in \R$ to $f$ -- and, needless to say, the larger $Q$, the tighter the lower bound. See Figure \ref{fig:lower_bound_comparison} for the results.

By Theorem~\ref{thm:properties_sls}, \ref{SLS} gives bounds that are at least as good as \ref{GLOB}. The plots show that SLS often gives strictly tighter bounds.

\subsubsection{Presentation of other bounds}
\label{subsubsec:other_bounds}

It is not straightforward to compare the performance of our lower bounds with bounds from the literature. In our setting, we compute a \textit{single} underestimator per instance -- which is then merely evaluated during the branch and bound process.\footnote{By fixing some variables at each node and then computing new underestimators, this could be improved but would need additional runtime for the computation of the new underestimator.} We could not find other underestimators with this property that give sensible results in branch and bound. However, there are lower bounds in the literature that are more general than ours since they consider restricted polynomial optimization problems and can hence be applied to any polynomial -- not only to those with positive definite leading form -- and are suitable for branch and bound if computed anew at each node. 
In addition to Algorithm~2 (with \ref{GLOB} and \ref{SLS}) we implemented the 
following four algorithms in a MATLAB framework for solving \ref{eq:problem}: 
three of them are branch and bound approaches as Algorithm~2 which use other bounds
(taken from \cite{buchheim2014box}, \cite{de2006integer}, and the continuous relaxation) while
our last algorithm is a simple brute force approach.

\begin{itemize}
	\item For arbitrary polynomials on boxes, Buchheim and D'Ambrosio~\cite{buchheim2014box} suggested to compute, for every term of $f$, the $L^1$-best separable underestimator. The sum of the underestimators is again separable, so its integer minimization is a univariate problem. For degree $d \leq 4$ and arbitrary $n$, they provide explicit underestimators. We hardcoded the explicit underestimators, and used the MATLAB builtins \verb+polyval+, \verb+polyder+ and \verb+roots+ to evaluate and differentiate the separable underestimators, and to compute their roots, respectively. As a suitable box at the subproblem $\mathcal{P}=(m, r_1, \ldots, r_m)$ we chose the box $[-L, L]^{n-m}$ where $L = \left \lfloor \sqrt[\leftroot{-2}\uproot{2}p]{R^p - |r_1|^p - \dots - |r_m|^{p}} \right \rfloor$. The authors suggest to successively halve the box into subboxes which does not fit into our scheme. This approach is abbreviated \textbf{SEP} in the plots.
	\item For nonnegative polynomials on polytopes $P$, De Loera et al.~\cite{de2006integer} approximate the maximum of $f$ on $P\cap \Z^n$ by the sequence $\sqrt[k]{\sum_{x \in P\cap \Z^n}f(x)^k}$. Each member of the sequence can be computed in polynomial time, using a reformulation as a limit of a rational function which in turn is based on the generating function of $P$. We did experiments with $k=2$ and $k=4$, the latter taking significantly longer, without giving much better results, so we restricted ourselves to $k=2$. Note that the suggested implementation uses residue techniques, while we just use symbolic limit computations. On the other hand, we improved the bounds as follows: To make their approach applicable to not necessarily nonnegative polynomials, the authors suggest to add the sufficiently large constant
		$$c:= \| f\|_0 \| f \|_{\infty} M^d$$
		to obtain $\overline f = f+c$ nonnegative on $P$. Here, $M \geq 0$ is a bound on the polyhedron s.t. $|x_i| \leq M$ for all $x \in P$; for $f = \sum_\alpha a_\alpha X^\alpha$, we use the zero ``norm'' $\| f\|_0~:=~\#\{ \alpha \ | \ a_\alpha \neq 0 \}$ and the infinity norm $\| f \|_\infty~:=~\max_\alpha \{ | a_\alpha | \}$. However, the constant $c':=\sum_{j=0}^d\|f_j\|_1 M^j$ suffices to ensure that $f+c'$ is nonnegative on $P$. A short calculation shows that $c'\leq c$ if $M \geq 1$, and in dense instances one often has $c' \ll c$. As polyhedron we again chose the box $[-L, L]^{n-m}$ from the previous bound. This bound is abbreviated to \textbf{POLYFIX} in the plots.
	\item We compute an sos approximation of the global continuous relaxation (\textbf{CR} in the plots) at each subproblem $\mathcal P = (m, r_1, \ldots, r_m)$, that is
		\begin{align*}
		\max \quad & \lambda \\
		\text{s.t.} \quad & f(r_1, \ldots, r_m, X_{m+1}, \ldots, X_{n}) - \lambda \text{ is sos in } \R[X_{m+1}, \ldots, X_{n}]
		\end{align*}
	\item Brute force enumeration with no lower bounds, abbreviated \textbf{BF}. As $f$ has to be evaluated at each node, we use \verb+matlabFunction+ to convert the Symbolic Math Toolbox object that encodes $f$ into a function handle that can be evaluated significantly faster.

	\item Algorithm~2 using \textbf{\ref{GLOB}} with parameters as described in Remark \ref{rem:input_to_glob_and_sls}.

	\item Algorithm~2 using \textbf{\ref{SLS}} with parameters as described in Remark \ref{rem:input_to_glob_and_sls}.
\end{itemize}

\subsubsection{Runtime comparison}

The implementation of the six different algorithms from Section~\ref{subsubsec:other_bounds} into our B\&B-framework gave the runtimes in Figure \ref{fig:runtimes} (logarithmic scale). On every instance each of the lower bounds had a maximum of 5 minutes to complete; if this time constraint was violated, the process was interrupted and the lower bound considered as unsuccessful on this instance. If the parameter $h$ could not be found by SOSTOOLS' \verb+findbound.m+ function, \ref{GLOB} and \ref{SLS} were considered to have violated the time constraint.

\begin{figure}[h!]
	\centering
	\subfloat{ \input{plotting/uest_experiment/uest_runtime_n_2-d_4.tex} } \qquad
	\subfloat{ \input{plotting/uest_experiment/uest_runtime_n_2-d_6.tex} } \\
	\vspace{0.5cm}
	\subfloat{ \input{plotting/uest_experiment/uest_runtime_n_3-d_4.tex} } \qquad
	\subfloat{ \input{plotting/uest_experiment/uest_runtime_n_4-d_2.tex} } \vspace{0.5cm} 
	\caption{Runtimes in $[\si{\second}]$.}
	\label{fig:runtimes}
\end{figure}
\begin{figure}[!h]
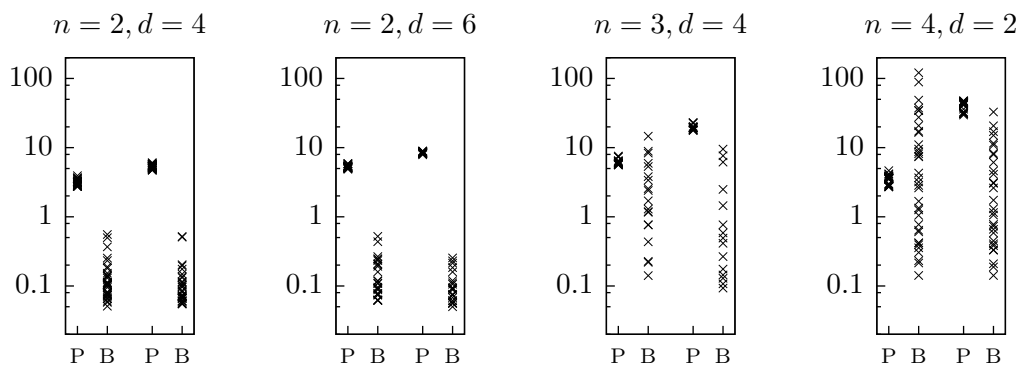

	\centering
	\subfloat{ \input{plotting/comparison_of_glob_and_sls/running_times-n_2-d_4.tex} } 
	\subfloat{ \input{plotting/comparison_of_glob_and_sls/running_times-n_2-d_6.tex} }
	\subfloat{ \input{plotting/comparison_of_glob_and_sls/running_times-n_3-d_4.tex} }
	\subfloat{ \input{plotting/comparison_of_glob_and_sls/running_times-n_4-d_2.tex} } \\
	\caption{Preprocessing (P) and B\&B (B) times -- \ref{GLOB} on the left, \ref{SLS} on the right.}
	\label{fig:prep_and_bandb_times}
\end{figure}

We infer from the plots that for a small number of variables, the problem size (i.e., $R$), is mostly so small that brute force is often the fastest approach. However, if instances get larger, brute force fails necessarily as the processing time is linear in the number of nodes. SEP is quite fast in small instances, but for large instances the running time deteriorates as an underestimator is computed at each node.
In our setting, POLYFIX takes too long to be competitive. The continuous relaxation is satisfactory for smaller instances but fails in some large instances. Concerning our bounds, in the two plots of Figure \ref{fig:runtimes} with $n=2$, there is a surprisingly little variance in runtime for GLOB and SLS. This can be explained from a further plot, see Figure~\ref{fig:prep_and_bandb_times}, in which we break down the preprocessing time, i.e., the time needed to compute a approximate continuous minimizer $h$ and the underestimator $g$, and the time needed for the actual branch and bound. It can be seen that the preprocessing time is more or less independent from the instance and takes in most instances significantly longer than the actual branch and bound. Also, it seems at first that \ref{SLS} takes mostly longer than \ref{GLOB}. However, this holds only true for the preprocessing phase: The corresponding sos program is larger, and so are preprocessing times. Indeed,  Figure \ref{fig:prep_and_bandb_times} reveals that \ref{GLOB} has shorter preprocessing times throughout, but is inferior in B\&B, as expected.

\section{Conclusion and Outlook}
\label{sec:ausblick}

In this paper we presented a new way of finding underestimators for
integer polynomial optimization and 
improved the bounds on the norm of integer and continuous minimizers. We implemented
both ideas within a branch \& bound approach showing how they improve
its performance. 
\medskip

Currently we compute one underestimator at the beginning of
the branch and bound process which is used for generating lower bounds
throughout the whole algorithm. Instead, one could also compute a new
underestimator at each node of the branch and bound tree. This would
improve the bounds but due to the comparatively large computation time for solving
an sos-program does not pay off in terms of overall efficiency. We currently
analyze in which nodes the computation of a new underestimators improves
the procedure. Along the lines of~\cite{Buch-Hueb-Sch13} 
we plan to analyze how to find an underestimator
which is likely to be a good one for all subnodes.

\medskip

We also point out that our procedure can be extended
to mixed-integer polynomial optimization: The norm bounds apply in the mixed-integer 
case as well, and we may use the proposed class of
underestimators, but with their mixed-integer minima (which are also
simple to obtain). Our underestimators can in principle also be used for 
\emph{constrained} polynomial optimization; however, 
it is subject to further investigation if the 
bounds provided are sharp enough for this case. 
We hence work on sos-programs which
provide underestimators which are able to take into account given constraints.

\appendix
\section{Computing the norm bounds}

In Remark \ref{rem:monoton} we saw that we get a tighter norm bound $R$ on the minimizers the closer the $c_j$ get to their optimal value $c_j^*$.
In the following, we present two means that improve on the approach \ref{item:cj_brute:1}. in Section \ref{subsec:the_bound} that do not rely on sos programming. The second method we present is a refinement of the first. For both, we improve the norm bound $R$ by replacing the estimate $|x^\alpha | \leq1$ on $\sphere^{n-1}_p$ with $|x^\alpha| \leq \hat x^\alpha$, where $\hat x$ is a continuous maximizer of the function $\sphere^{n-1}_p \rightarrow \R$, $x \mapsto x^\alpha$. 
\subsection{A direct improvement}

One has the following closed form for the continuous minimizer $\hat x$ with nonnegative coordinates:

\begin{lem}
	\label{lem:maximum_of_monomial_on_sphere}
	Let $0 \neq \alpha \in \N_0^n$ and $p \in [1, \infty)$. Then, the monomial $X^\alpha$ attains its maximum on $\sphere^{n-1}_p$ at $\hat x$ with coordinates
	\begin{equation}
		\hat x_i = \sqrt[\leftroot{-2}\uproot{2}p]{\frac{\alpha_i}{\sum_{i=1}^n \alpha_i}}, \quad i = 1, \ldots, n.
		\label{eq:maximum_on_sphere}
	\end{equation}
		\begin{proof}
By a simple analysis, the proof can be reduced to $\alpha_i \geq 1$ for $i=1, \ldots, n$ and then to maximization of $X^\alpha$ on $\{ x \in \sphere^{n-1}_p \ | \ x_1 > 0, \ldots, x_n > 0 \}$. Using the method of Lagrange multipliers, the claim follows from a short calculation.
		\end{proof}
\end{lem}

\begin{obs}
	\label{obs:cj_refined} Denote by $\hat x_{(\alpha)}$ the maximizer of $X^\alpha$ on $\sphere^{n-1}_p$ as in (\ref{eq:maximum_on_sphere}). Hence for $x \in \sphere^{n-1}_p$ we have
		\begin{equation}
			f_j(x) = \sum_{|\alpha| = j} a_\alpha x^{\alpha} \geq \sum_{|\alpha| = j} - |a_\alpha| \cdot (\hat x_{(\alpha)})^\alpha =: c_j.
			\label{eq:cj_lagrange}
		\end{equation}
		This $c_j$ is as least as large as approach \ref{item:cj_brute:1}. from Section \ref{subsec:the_bound} since, for $0 \neq \alpha$, $(\hat x_{(\alpha)})^\alpha < 1$ -- unless $X^\alpha \in \R[X_i]$ for some $i$, in which case $\hat x_{(\alpha)} = e_i$, the $i$-th unit vector, and thus $(\hat x_{(\alpha)})^\alpha = 1$.
\end{obs}

\subsection{A different approach}

This last approach on computing bounds $c_j$ is different to the ones before, as we actually compute $2^n$ norm bounds: We restrict $f$ to each of the $2^n$ orthants
$$H_\tau = \{x \in \R^n \ |\ \tau_i x_i \geq 0\} \text{ for }\tau \in \{-1,1\}^n$$ 
and compute norm bound on integer minimizers of every $f|_{H_\tau}$.
This has the advantage that, roughly speaking, we may neglect half of the terms of $f = \sum a_\alpha X^\alpha$. Also, minimization on $H_\tau$ can be reduced to minimization on $H_{(1, \ldots, 1)}$, i.e., the set of those $x \in \R^n$ with $x \geq 0$, as we shall see in a moment. 

Introducing the notation $|a|^- = |\min(a, 0)|$ for $a \in \R$ and with $\hat x$ from \eqref{eq:maximum_on_sphere}, we have for every term
$a_\alpha x^\alpha \geq - |a_\alpha|^- x^\alpha \geq -| a_\alpha|^- \hat x^\alpha$ as $x \geq 0$, thus
$$f_j(x) = \sum_{|\alpha| = j } a_\alpha x^\alpha \geq \underbrace{\sum_{|\alpha| = j} -|a_\alpha|^- \hat x^\alpha}_{=:c_j^{(1, \ldots, 1)}}, \quad x \in \sphere^{n-1}_p \text{ and } x \geq 0,$$
which means about half of the coefficients are neglected in comparison to \eqref{eq:cj_lagrange}, if signs are distributed equally among the $a_\alpha$. Now let $R^{(1, \ldots, 1)}$ be the largest real root of
$$q^{(1, \ldots, 1)}(\lambda) := c_d \lambda^d + \sum_{j=1}^{d-1} c_j^{(1, \ldots, 1)} \lambda^j.$$ 
The verbatim argument of Theorem \ref{thm:bound_on_norm} shows that $f(x) > f(0)$ for $\| x \|_p > R^{(1, \ldots, 1)}$ and $x \geq 0$. This bounds integer and continuous minimizers on $H_{(1, \ldots, 1)}$.
Bounding the norm of minimizers of $f$ on $H_\tau$, $\tau \in \{-1,1\}^{n}$, can be reduced to bounding the norm of minimizers on $H_{(1, \ldots, 1)}$ by a simple change of coordinates. To this end, let $\tau(x) = (\tau_1 x_1, \ldots, \tau_n x_n)$, $x \in \R^n$, and $f^\tau$ be the polynomial
$$f^{\tau}(x) := f(\tau(x)) = \sum_{\alpha }a_\alpha \tau^\alpha x^\alpha, \quad \tau \in \{-1, 1 \}^n.$$

As $\tau^\alpha \in \{-1, 1\}$, $f$ and $f^\tau$ merely differ in the sign of their coefficients, and $f^\tau_d(x) \geq c_d$ still holds for $x \in \sphere^{n-1}_p$ as the sphere is $\tau$-invariant, that is $\tau(\sphere^{n-1}_p) = \sphere^{n-1}_p$. Similarly to before, denote by $R^{\tau}$ the largest real root of
$$q^{\tau}(\lambda) = c_d \lambda^d + \sum_{j=1}^{d-1} c_j^{\tau} \lambda^{j},$$
with $c_j^\tau = -|a_\alpha \tau^{\alpha}|^-  \hat x^\alpha$.
It is now clear that $f^{\tau}(x) > f(0)$ for $\| x\|_p > R^\tau$ and $x \geq 0$, equivalently, $f(x) > f(0)$ for $\| x\|_p > R^\tau$ and $x \in H_\tau$.

\medskip
This results in more effort in the preprocessing, but reduces the number of feasible solutions.

\bibliography{polynom,../eigen,../t}

\end{document}